\def\ver{Nov. 18, 2002, v.14}
\documentstyle{amsppt}
\magnification=1200
\hoffset=-10pt
\voffset=-20pt
\hsize=6.5truein
\vsize=8.9truein
\topmatter
\title Bloch's Conjecture, Deligne Cohomology\\
and Higher Chow Groups
\endtitle
\author Morihiko Saito
\endauthor
\affil RIMS Kyoto University, Kyoto 606-8502 Japan \endaffil
\keywords Deligne cohomology, higher Chow group,
cycle map
\endkeywords
\subjclass 14C30, 32S35\endsubjclass
\abstract We express the kernel of Griffiths' Abel-Jacobi
map by using the inductive limit of Deligne cohomology in
the generalized sense (i.e. the absolute Hodge cohomology of
A. Beilinson). This generalizes a result of L. Barbieri-Viale
and V. Srinivas in the surface case. We then show that the
Abel-Jacobi map for codimension 2 cycles and the Albanese
map are bijective if a general hyperplane section
is a surface for which Bloch's conjecture is proved.
In certain cases we verify Nori's conjecture on the Griffiths
group. We also prove a weak Lefschetz-type theorem for (higher)
Chow groups, generalize a formula for the Abel-Jacobi map
of higher cycles due to Beilinson and Levine to the smooth
non proper case, and give a sufficient condition for the
nonvanishing of the transcendental part of the image by the
Abel-Jacobi map of a higher cycle on an elliptic surface,
together with some examples.
\endabstract
\endtopmatter
\tolerance=1000
\baselineskip=12pt
\def\scirc{\raise.2ex\hbox{${\scriptstyle\circ}$}}
\def\ssbull{\raise.2ex\hbox{${\scriptscriptstyle\bullet}$}}
\def\cssbull{\,\raise.4ex\hbox{${\scriptscriptstyle\bullet}$}\,}
\def\moplus{\raise.1ex\hbox{$\bigoplus$}}
\def\msum{\raise.1ex\hbox{$\sum$}}
\def\anD{{\langle D\rangle}}
\def\bA{{\Bbb A}}
\def\bC{{\Bbb C}}
\def\bD{{\Bbb D}}
\def\bG{{\Bbb G}}
\def\bP{{\Bbb P}}
\def\bQ{{\Bbb Q}}
\def\bR{{\Bbb R}}
\def\bZ{{\Bbb Z}}
\def\cC{{\Cal C}}
\def\cD{{\Cal D}}
\def\cE{{\Cal E}}
\def\cF{{\Cal F}}
\def\cH{{\Cal H}}
\def\cK{{\Cal K}}
\def\cO{{\Cal O}}
\def\cU{{\Cal U}}
\def\cV{{\Cal V}}
\def\cZ{{\Cal Z}}
\def\og{\overline{g}}
\def\oi{\overline{i}}
\def\oX{\overline{X}}
\def\oZ{\overline{Z}}
\def\oga{\overline{\gamma}}
\def\oze{\overline{\zeta}}
\def\tzero{\widetilde{0}}
\def\tD{\widetilde{D}}
\def\tf{\widetilde{f}}
\def\tg{\widetilde{g}}
\def\tX{\widetilde{X}}
\def\tY{\widetilde{Y}}
\def\tZ{\widetilde{Z}}
\def\tga{\widetilde{\gamma}}
\def\Griff{\hbox{{\rm Griff}}}
\def\Zar{\text{{\rm Zar}}}
\def\Rat{\hbox{{\rm Rat}}}
\def\Dec{\hbox{{\rm Dec}}\,}
\def\Pic{\hbox{{\rm Pic}}}
\def\rat{\hbox{{\rm rat}}}

\def\div{\hbox{{\rm div}}\,}
\def\codim{\hbox{{\rm codim}}}
\def\Tr{\hbox{{\rm Tr}}}
\def\NS{\hbox{{\rm NS}}}
\def\IC{\hbox{{\rm IC}}}
\def\CH{\hbox{{\rm CH}}}
\def\Spec{\hbox{{\rm Spec}}}
\def\Ext{\hbox{{\rm Ext}}}
\def\Gr{\text{{\rm Gr}}}
\def\Re{\hbox{{\rm Re}}}
\def\Im{\hbox{{\rm Im}}}
\def\Ker{\hbox{{\rm Ker}}}
\def\Coker{\hbox{{\rm Coker}}}
\def\Hom{\hbox{{\rm Hom}}}
\def\Alb{\hbox{{\rm Alb}}}
\def\MHM{\text{{\rm MHM}}}
\def\MHS{\text{{\rm MHS}}}

\def\Sing{\hbox{{\rm Sing}}\,}
\def\supp{\hbox{{\rm supp}}\,}
\def\Hv{\text{\rm Hv}}
\def\an{\text{\rm an}}
\def\BM{\text{\rm BM}}
\def\rel{\text{\rm rel}}
\def\free{\text{\rm free}}
\def\ab{\text{\rm ab}}
\def\AJ{\text{\rm AJ}}
\def\alg{\text{\rm alg}}
\def\ind{\text{\rm ind}}
\def\dec{\text{\rm dec}}
\def\tor{\text{\rm tor}}
\def\reg{\text{\rm reg}}
\def\hom{\text{\rm hom}}
\def\id{\text{\rm id}}
\def\simto{\buildrel\sim\over\to}
\def\SameAuthor{\vrule height3pt depth-2.5pt width1cm}

\document
\centerline{{\bf Introduction}}
\bigskip
\noindent
Let
$ X $ be a smooth projective complex surface.
D. Mumford [41] showed that the kernel of the Albanese map
$ \CH_{0}(X)^{0} \to \Alb(X) $ is `huge' if
$ X $ has a nontrivial global
$ 2 $-from (i.e. if
$ p_{g}(X) \ne 0 $).
Then S. Bloch [9] conversely conjectured

\medskip\noindent
(0.1)\,\,\,
The Albanese map
$ \CH_{0}(X)^{0} \to \Alb(X) $ is injective if
$ p_{g}(X) = 0 $.

\medskip\noindent
This conjecture was proved in [13] if X is not of general type,
but the general case still remains open, see [3], [31], [54], etc.
Related to this, L. Barbieri-Viale and V. Srinivas [2]
(see also [29], [46]) constructed an exact sequence
$$
H_\cD^{3}(X,\bZ(2)) \to \underset{U}\to\varinjlim\,
H_\cD^{3}(U,\bZ(2)) \to \CH_{0}(X)^{0} \to \Alb(X),
$$
where
$ U $ runs over the nonempty open subvarieties of
$ X $, and
$ H_\cD^{3}(U,\bZ(2)) $ denotes Deligne cohomology.
This describes the kernel of the Albanese map,
and follows also from the local-to-global spectral sequence in the
theory of Bloch-Ogus [14].

In this paper, we generalize this to the higher dimensional case,
using Deligne cohomology in a generalized sense.
The notion of Deligne cohomology was first introduced by Deligne in the
case
$ X $ is smooth and proper.
It is a natural generalization of the first two terms of the exponential
sequence (see [25]).
The generalization to the open or singular case
was first done by A. Beilinson [4] and
H. Gillet [29], where the weight filtration was not used.
Later Beilinson [5] found a more natural generalization from the view
point of mixed Hodge theory, which he calls absolute Hodge cohomology,
and denotes by
$ H_\cH^{i}(X,A(k)) $,
$ H_{\cH^{p}}^{i}(X,A(k)) $,
where
$ A $ is a subring of
$ \bR $.
In this paper we denote them by
$ H_\cD^{i}(X,A(k))' $,
$ H_\cD^{i}(X,A(k))'' $ respectively, see (1.1).

Let
$ \CH_{\hom}^{p}(X) $ be the subgroup of
$ \CH^{p}(X) $ consisting of cycles homologically equivalent to zero.
There is Griffiths' Abel-Jacobi map to the intermediate Jacobian
$ J^{p}(X) $.
Its kernel is described by using Deligne cohomology as follows:

\medskip\noindent
{\bf 0.2.~Theorem.} {\it
Let
$ X $ be a smooth proper variety over
$ \bC $.
For any integer
$ p $,
there is a canonical exact sequence
$$
H_\cD^{2p-1}(X,\bZ(p)) \to \underset{Y}\to\varinjlim\,
H_\cD^{2p-1}(X \setminus Y,\bZ(p))''
\to \CH_{\hom}^{p}(X) \to J^{p}(X),
$$
where the inductive limit is taken over the closed subvarieties
$ Y $ of
$ X $ with pure codimension
$ p - 1 $.
Here
$ H_\cD^{2p-1}(X \setminus Y,\bZ(p))'' $ may be replaced by
$ H_\cD^{2p-1}(X \setminus Y,\bZ(p))' $ in general, and by
$ H_\cD^{2p-1}(X \setminus Y,\bZ(p)) $ if
$ p = \dim X $.
}

\medskip
For the proof of this, we study the cycle map to Deligne homology
of a singular variety.
Using the compatibility of the cycle map with the localization
sequence (2.7), Theorem (0.2) is reduced to

\medskip\noindent
{\bf 0.3.~Theorem.} {\it
For a variety
$ Y $ of pure dimension
$ m $, the cycle map induces isomorphisms
$$
cl : \CH^{1}(Y) \simto
H_{2m-2}^\cD(Y,\bZ(m-1))'' =
H_{2m-2}^\cD(Y,\bZ(m-1))'.
$$
If
$ m = 1 $, this holds also for
$ H_{2m-2}^\cD(Y,\bZ(m-1)) $.
}

\medskip
Indeed, using the localization sequence,
we can reduce it to the smooth case (see (3.1) below),
because a similar isomorphism for
$ \CH^{1}(Y,1) $ is already known, see [32, 3.1].
It is also possible to describe
$ \CH_{\hom}^{1}(Y) $ by using the normalization of
$ Y $, see (3.5).
This may be useful for explicit calculation.

As another application of (0.3), we prove (see (3.13)):

\medskip\noindent
{\bf 0.4.~Theorem.} {\it
Let
$ X $ be a smooth projective complex algebraic variety of
dimension
$ n $, and
$ Y $ be an intersection of
$ n - 2 $ general hyperplane sections of
$ X $.
Assume
$ Y $ has no global
$ 2 $-forms and the Albanese map is injective for
$ Y $, i.e. Bloch's conjecture {\rm (0.1)} holds.
Then the Albanese map for
$ X $ and the Abel-Jacobi map for cycles of codimension
$ 2 $ on
$ X $ are bijective.
}

\medskip
The proof of the injectivity of the Albanese map is
relatively easy,
and it implies the surjectivity of the Abel-Jacobi map for
cycles of codimension
$ 2 $, see [26].
For the injectivity of the latter we show (see (3.11)):

\medskip\noindent
{\bf 0.5.~Theorem.} {\it
Let
$ X $ be an irreducible smooth proper complex algebraic variety.
Let
$ f : X \to S $ be a surjective morphism to a smooth variety
$ S $.
Assume that general fibers
$ X_{s} $ of
$ f $ are connected and have no global
$ 2 $-forms, and the Abel-Jacobi map for cycles of codimension
$ 2 $ on general fibers is injective.
If
$ H^{1}(X_{s},\bQ) \ne 0 $ for a general
$ s \in S $, we assume further that
$ S = \bP^{1} $ and the restriction morphism
$ H^{1}(X,\bQ) \to H^{1}(X_{s},\bQ) $ is an isomorphism for a
general
$ s \in S $.
Then the kernel of the Abel-Jacobi map tensored with
$ \bQ $ for cycles of codimension
$ 2 $ on
$ X $ comes from that on
$ S $.
}

\medskip
This implies under the assumption of (0.5) that Nori's conjecture
[42] on the Griffiths group [30] for
$ X $ is true if it holds for
$ S $ (e.g. if
$ \dim S \le 2 $), see also (3.14) below.
Recall that the conjecture predicts an isomorphism between the
Griffiths group for cycles of codimension
$ 2 $ and the quotient of the image of the Abel-Jacobi map
divided by the maximal abelian subvariety (and is equivalent
to that Abel-Jacobi equivalence is stronger than algebraic
equivalence).
This conjecture can be deduced from a well-known conjecture of
Beilinson [6] and Bloch [9] on a conjectural filtration of the
Chow groups (assuming the Hodge conjecture).
The hypothesis on the vanishing of
$ H^{1}(X_{s},\bQ) $ in (0.5) is satisfied if general fibers
$ X_{s} $ are surfaces of general type (see e.g. [52]).

Under the assumption of (0.4),
$ X $ has no global
$ p $-forms for
$ p > 1 $,
and this is compatible with [44].
It is conjectured that the injectivity of the Abel-Jacobi map
for cycles of codimension
$ 2 $ on
$ X $ in (0.4) should hold by assuming only that
$ X $ has no global
$ 2 $-form.
The hypothesis on Bloch's conjecture (0.1) in (0.4) is satisfied
at least if
$ Y $ is not of general type,
see [13] (and also [3], [31], [54], etc.)
The assumption on
$ X $ in (0.4) is satisfied for example by cubic threefolds
(see also [17]) and smooth complete intersections of degree
$ (2,2) $ in
$ \bP^{5} $.
There are other examples because any smooth projective variety
is a general hyperplane section of any
$ \bP^{1} $-bundle having a section over it
(choosing a projective embedding appropriately).
In some cases we can show that algebraic and homological
equivalences coincide for cycles of codimension
$ 2 $, see [1] and also (3.17) below.
Note that the Griffiths group is not finitely generated for
general Calabi-Yau threefolds [56].

As another application we have a weak Lefschetz-type theorem
for Chow groups (see (3.15), and for a similar assertion about
higher cycles, see (3.16)):

\medskip\noindent
{\bf 0.6.~Theorem.} {\it
Let
$ X $ be a smooth projective complex algebraic variety.
Take a Lefschetz pencil
$ f : \tX \to S := \bP^{1} $ where
$ \pi : \tX \to X $ is the blow-up along an
intersection of two generic hyperplane sections.
Let
$ S' $ be any nonempty open subvariety of
$ S $ over which
$ f $ is smooth.
Assume
$ \dim X \ge 4 $.
Then
$ \zeta \in \CH_{\hom}^{2}(X) $ is zero if its restriction to
$ X_{s} := f^{-1}(s) $ vanishes for any
$ s \in S' $.
}

\medskip
Concerning Bloch's conjecture (0.1), it is known [33] that the
conjecture is related to the surjectivity of the cycle map
$$
cl : \CH^{p}(X,1) \to H_{\cD}^{2p-1}(X,\bZ(p))''
$$
for certain smooth nonproper varieties
$ X $ (see also (3.3) below).
Here
$ \CH^{p}(X,m) $ denotes Bloch's higher Chow group [10].
In the smooth proper case, Beilinson [5] and Levine [37] described
the cycle map explicitly by using currents like for Griffiths'
Abel-Jacobi map.
In this paper, we extend this to an explicit description in the smooth
nonproper case, see (4.3).

Let us return to the case of a smooth proper surface
$ X $.
By [15], the conclusion of (0.1) would imply the decomposability of
$ \CH^{2}(X,1)_{\bQ} $ (i.e. it is generated by the image of
$ \Pic(X)\otimes_{\bZ}\bC^{*}) $, see also [26].
So it is conjectured that
$ \CH^{2}(X,1)_{\bQ} $ is decomposable if
$ p_{g}(X) = 0 $.
Thus it would be interesting whether the reduced higher Abel-Jacobi map
$$
\CH_{\ind}^{2}(X,1)_{\bQ} \to
J(H^{2}(X,\bZ)(2))_{\bQ}/\NS(X)_{\bQ}\otimes_{\bZ}\bC^{*},
\leqno(0.7)
$$
which is induced by the above cycle map, is injective in general.
Here
$ \CH_{\ind}^{2}(X,1)_{\bQ} $ is the quotient of
$ \CH^{2}(X,1)_{\bQ} $ by the image of
$ \Pic(X)_{\bQ}\otimes_{\bZ}\bC^{*} $, see (2.2.3) and (4.5.1) below.

This injectivity is related to Voisin's conjecture [55] on the
countability of
$ \CH_{\ind}^{2}(X,1)_{\bQ} $, because the image of the reduced
higher Abel-Jacobi map (0.7) is countable [40].
The kernel of this map is isomorphic to
$$
\Coker(K_{2}(\bC(X))_{\bQ} \to
\underset{U}\to\varinjlim\, \Hom_{\MHS}(\bQ, H^{2}(U,\bQ)(2)))
$$
by [43], [46], where the morphism of
$ K_{2}(\bC(X))_{\bQ} $ is given by
$ d \log \wedge d \log $ at the level of integral logarithmic forms,
and the inductive limit is taken over the nonempty open subvarieties
$ U $ of $ X $ (see also [5], 6.1).
This isomorphism follows easily from the localization sequence of
mixed Hodge structures together with the fact that the residue of
$ d \log f \wedge d \log g $ coincides with the logarithmic
differential of the tame symbol of
$ \{f,g\} $ up to sign.
It holds also for open subvarieties
$ U $ if
$ H^{3}(U,\bQ) = 0 $.

In view of these considerations we are interested in constructing
examples of indecomposable higher cycle such that the transcendental
part of its image by the higher Abel-Jacobi map does not vanish
(i.e. its image is not contained in the image of
$ F^{1}H^{2}(X,\bC) $ in the Jacobian).
We give a sufficient condition for it together with
some examples in the case of elliptic surfaces, see (5.2--3) below.
(In an earlier version of this paper, such an example was constructed
by calculating period integrals of elliptic curves and using
double integration.)
Another example satisfying the above property is found independently by
P. del Angel and S. M\"uller-Stach [21].
The support of our cycle is irreducible, and is a fiber of an elliptic
surface.
Such an example does not seem to have appeared in the literature.

I would like to thank A. Rosenschon and L. Barbieri-Viale for useful
discussions.

\medskip
In Sect.~1, we review some basic facts from
the theories of Deligne cohomology and mixed Hodge Modules
which are needed in this paper.
In Sect.~2, we recall the definition of Bloch's higher Chow groups
and the cycle map.
In Sect.~3 we prove the main theorems (0.2--5) and some related
assertions.
In Sect.~4 we describe explicitly the cycle map for higher cycles,
and construct examples of indecomposable higher cycles in
Sect. 5.

In this paper a variety means a separated scheme of finite type over
$ \bC $.
All sheaves are considered on the associated analytic spaces, and
$ H^{j}(X^{\an},\bQ) $ is denoted by
$ H^{j}(X,\bQ) $.

\bigskip\bigskip
\centerline{{\bf 1. Deligne cohomology and mixed Hodge Modules}}

\bigskip
\noindent
{\bf 1.1. ~Deligne cohomology} (see [4, 5, 21, 23, 25, 27, 32], etc.)
Let
$ X $ be a smooth variety, and
$ \oX $ a smooth compactification of
$ X $ such that
$ D := \oX \setminus X $ is a divisor with normal crossings.
Let
$ j : X \to \oX $ denote the inclusion morphism.
Let
$ A $ be
$ \bZ $ or
$ \bQ $ for simplicity in this paper.
Then
$ A $-Deligne cohomology
$ H_{\cD}^{i}(X,A(k)) $ is defined to be the
$ i $-th hypercohomology group of
$$
C_{\oX\anD}^{\ssbull}A(k) := C(\bold{R}j_{*}
A_{X}(k)\oplus \sigma_{\ge k}{\Omega}_{\oX}^{\ssbull}(\log D) \to
\bold{R}j_{*}{\Omega}_{X}^{\ssbull})[-1],
$$
where
$ A_{X}(k) = (2\pi i)^{k}A_{X} \subset \bC_{X} $.

For a complex variety
$ X $ in general, let
$ K= (K_A,(K_\bQ,W),(K_\bC,F,W)) $ be the complex of
graded-polarizable mixed
$ A $-Hodge structures corresponding by [5, 3.11] to the mixed
Hodge complex calculating the cohomology of
$ X $ which is defined by using a simplicial resolution of a
compactification of
$ X $ as in [22].
Let
$ \MHS(A)^{p} $ (resp.
$ \MHS(A) $) denote the abelian category of graded-polarizable
(resp. not necessarily graded-polarizable) mixed
$ A $-Hodge structures.
Then we define
$ A $-Deligne cohomology in the generalized sense by
$$
\aligned
H^{i}_\cD(X,A(k))
&= H^{i}(C(K_A(k)\oplus
F^{k}K_\bC \to K_\bC)[-1]),
\\
H^{i}_\cD(X,A(k))'
&= \bold{R}\Hom_{\MHS(A)}(A,K(k)[i]),
\\
H^{i}_\cD(X,A(k))''
&= \bold{R}\Hom_{\MHS(A)^{p}}(A,K(k)[i]).
\endaligned
$$
The last two are called absolute Hodge cohomology in [5],
and denoted respectively by
$ H_\cH^{i}(X,A(k)) $,
$ H_{\cH^{p}}^{i}(X,A(k)) $.
By loc. cit. we have natural morphisms
$$
H^{i}_\cD(X,A(k))'' \to
H^{i}_\cD(X,A(k))' \to H^{i}_\cD(X,A(k)).
$$

Similarly, let
$ K'= (K'_A,(K'_\bQ,W),(K'_\bC,F,W)) $ be the dual of the
complex of graded-polarizable mixed
$ A $-Hodge structures corresponding by [5, 3.11] to the mixed
Hodge complex calculating the cohomology with compact support of
$ X $ which is defined by using a simplicial resolution of a
compactification of
$ X $ together with that of the divisor at infinity.
Then we define
$ A $-Deligne homology in the generalized sense by
$$
\aligned
H_{i}^\cD(X,A(k))
&= H^{-i}(C(K'_A(-k)\oplus
F^{-k}K'_\bC \to K'_\bC)[-1]),
\\
H_{i}^\cD(X,A(k))'
&= \bold{R}\Hom_{\MHS(A)}(A,K'(-k)[-i]),
\\
H_{i}^\cD(X,A(k))''
&= \bold{R}\Hom_{\MHS(A)^{p}}(A,K'(-k)[-i]).
\endaligned
$$
We have also natural morphisms
$$
H_{i}^\cD(X,A(k))'' \to
H_{i}^\cD(X,A(k))' \to H_{i}^\cD(X,A(k)).
$$

If
$ X $ is smooth of pure dimension
$ n $, then
$ K = K'(n)[2n] $ so that
$$
H_{i}^{\cD}(X,A(k)) = H_{\cD}^{2n-i}(X,A(n-k)),
\leqno(1.1.1)
$$
and similarly for
$ H_{i}^{\cD}(X,A(k))' $,
$ H_{i}^{\cD}(X,A(k))'' $.

For a mixed
$ A $-Hodge structure
$ H = (H_A, (H_\bQ,W), (H_\bC,F,W)) $ and an integer
$ k $, we define
$$
\aligned
J(H(k))
&= H_\bC/(H_A(k)_{\free} + F^{k}H_\bC),
\\
J'(H(k))
&= W_{2k}H_\bC/((W_{2k}H_A)(k)_{\free} + F^{k}H_\bC),
\\
J''(H(k))
&= W_{2k-1}H_\bC/(((W_{2k}H_A)(k)_{\free} + F^{k}H_\bC)
\cap W_{2k-1}H_\bC),
\endaligned
$$
where
$ H_A(k)_{\free} = H_A(k)/H_A(k)_{\tor} $.
We define also
$$
F^{k}W_{2k}H_{A}(k) = \Ker(H_{A}(k) \to
H_\bC/F^{k}W_{2k}H_\bC).
$$
(Similarly for
$ F^{k}H_{A}(k) $.)
Then we have short exact sequences
$$
\aligned
0 \to J(H^{i-1}(X,A)(k))
&\to H_{\cD}^{i}(X,A(k)) \to
F^{k}H^{i}(X,A)(k) \to 0,
\\
0 \to J(H^{i-1}(X,A)(k))'
&\to H_{\cD}^{i}(X,A(k))' \to
F^{k}W_{2k}H^{i}(X,A)(k) \to 0,
\\
0 \to J(H^{i-1}(X,A)(k))''
&\to H_{\cD}^{i}(X,A(k))'' \to
F^{k}W_{2k}H^{i}(X,A)(k) \to 0,
\endaligned
\leqno(1.1.2)
$$
because
$ J'(H(k)) = \Ext_{\MHS(A)}^{1}(A,H(k)) $,
$ J''(H(k)) = \Ext_{\MHS(A)^{p}}^{1}(A,H(k)) $ by [16]
and the semisimplicity of polarizable Hodge structures.
We have also
$$
\aligned
0 \to J(H_{i+1}^{\BM}(X,A)(-k))^{\prime\prime}
&\to H^{\cD}_{i}(X,A(k))^{\prime\prime}
\\
&\quad\to F^{-k}W_{-2k}H_{i}^{\BM}(X,A)(-k) \to 0,
\endaligned
\leqno(1.1.3)
$$
etc.
Here
$ H_{i}^{\BM}(X,A) $ denotes Borel-Moore homology.

It is known that
$ H_{\cD}^{i}(X,A(k)) $ and
$ H_{i}^{\cD}(X,A(k)) $ (together with Deligne local
cohomology) satisfy the axioms of Bloch-Ogus [14],
see [4], [29], [32], etc.
In particular, we have a canonical long exact sequence
$$
\to H_{i}^{\cD}(Y,A(k))
\to H_{i}^{\cD}(X,A(k))
\to H_{i}^{\cD}(U,A(k))
\to H_{i-1}^{\cD}(Y,A(k)) \to
\leqno(1.1.4)
$$
for a closed subvariety
$ Y $ of
$ X $ and
$ U = X \setminus Y $.
(This is functorial for
$ Y $,
$ U $.)
Similar assertions hold for
$ H_{i}^{\cD}(X,A(k))' $,
$ H_{i}^{\cD}(X,A(k))'' $.
(The assertion for
$ H_{i}^{\cD}(X,\bQ(k))'' $ follows also from [49]
using (1.6) below.)
We can similarly define Deligne local cohomology supported on
$ Y \subset X $.
If
$ X $ is smooth this coincides with Deligne homology of
$ Y $.

\medskip\noindent
{\bf 1.2.~Remark.} Let
$ K = (K_{\bZ}, (K_{\bQ},W), (K_{\bC};F,W);
K'_{\bQ}, (K'_{\bC},W)) $ be a polarizable
mixed Hodge complex in the sense of [5, 3.9].
They are endowed with (filtered) quasi-isomorphisms
$$
\aligned
\alpha_{1} : K_{\bZ}\otimes_{\bZ}\bQ \to
K'_{\bQ},\quad
&\alpha_{2} : K_{\bQ} \to K'_{\bQ},
\\
\alpha_{3} : (K_{\bQ},W)\otimes_{\bQ}\bC \to
(K'_{\bC},W),\quad
&\alpha_{4} : (K_{\bC},W) \to (K'_{\bC},W)
\endaligned
$$
such that
$ (\Gr_{i}^{W}K_{\bQ}, \Gr_{i}^{W}(K_{\bC},F)) $ together with the
isomorphism in the derived category
$ \Gr_{i}^{W}K_{\bQ}\otimes_{\bQ}\bC =
\Gr_{i}^{W}K_{\bC} $ induced by
$ \alpha_{3}, \alpha_{4} $ is a polarizable Hodge complex of weight
$ i $ in the sense of [22], i.e.
$ \Gr_{i}^{W}(K_{\bC},F) $ is strict and
$ H^{j}(\Gr_{i}^{W}K_{\bQ}, \Gr_{i}^{W}(K_{\bC},F)) $ is a
polarizable Hodge structure of weight
$ i + j $.

Let
$ \Dec W $ be as in loc. cit.
By definition, we have a canonical surjection
$$
(\Dec W)_{0}K_{\bQ}^{j} \to H^{j}\Gr_{-j}^{W}K_{\bQ},
$$
(and similarly for
$ K_{\bC}, K'_{\bC}) $.
For a
$ \bQ $-Hodge structure
$ H = (H_{\bQ}, (H_{\bC},F)) $ of weight
$ 0 $,
let
$$
H^{(0)} = \Hom_{\MHS}(\bQ,H),
$$
which is identified with a subgroup of
$ H_{\bQ}, H_{\bC} $.
We define the subcomplex
$ (\Dec W)_{0}^{(0)}K_{\bQ} $ of
$ (\Dec W)_{0}K_{\bQ} $ so that
$ (\Dec W)_{0}^{(0)}K_{\bQ}^{j} $ is the inverse image of
$ (H^{j}\Gr_{-j}^{W}K)^{(0)} $ by the above morphism (and similarly for
$ (\Dec W)_{0}^{(0)}K_{\bC} $,
$ (\Dec W)_{0}^{(0)}K'_{\bC}) $.

Let
$ (W_{0}H^{j}K_{\bQ})^{(0)} $ be the inverse image of
$ (\Gr_{0}^{W}H^{j}K_{\bQ})^{(0)} $ by the projection of
$ W_{0}H^{j}K_{\bQ} $ to
$ \Gr_{0}^{W}H^{j}K_{\bQ} $.
Since
$ d_{1} $ of the weight spectral sequence is a morphism of Hodge
structures, we can show the canonical quasi-isomorphism
$$
\tau_{\le j}(\Dec W)_{0}^{(0)}K_{\bQ}/
\tau_{<j}(\Dec W)_{0}^{(0)}K_{\bQ} \to
(W_{0}H^{j}K_{\bQ})^{(0)},
$$
and similarly for
$ K_{\bC}, K'_{\bC} $.

We define a complex
$ \Gamma (D''_{\cH}K) $ to be the single complex associated with
$$
K_{\bZ}\oplus (\Dec W)_{0}^{(0)}K_{\bQ}\oplus
F^{0}(\Dec W)_{0}^{(0)}K_{\bC}
\overset\phi\to\to
K'_{\bQ}\oplus (\Dec W)_{0}^{(0)}K'_{\bC},
$$
where
$ \phi $ is induced by
$ (\alpha_{1} - \alpha_{2})\oplus (\alpha_{3} - \alpha_{4}) $,
and the degree of the source of
$ \phi $ is zero.
Then, by an argument similar to [5], we can show the isomorphism
$$
\Hom_{\cD''}(\bZ,K) = H^{0}\Gamma (D''_{\cH}K),
\leqno(1.2.1)
$$
where
$ \cD'' $ denotes the category of polarizable mixed Hodge
complexes in the sense of [5, 3.9].
So we can define
$ H_{\cD}^{i}(X,\bZ(k))'' $ taking a mixed Hodge complex
which calculates the cohomology of
$ X $ as in [22].
Note that (1.2.1) implies the equivalence of categories
$$
D^{b}\MHS(\bZ)^{p} \simto \cD''
$$
in Lemma 3.11 of [5], and that it is easy to show the exact sequence
$$
0 \to \Ext_{\MHS(\bZ)^{p}}^{1}(\bZ, H^{-1}K)
\to \Hom_{\cD''}(\bZ,K) \to
\Hom_{\MHS(\bZ)}(\bZ, H^{0}K) \to 0,
$$
using the truncation
$ \tau $.

We can similarly define
$ \Gamma (D'_{\cH}K), \Gamma (D_{\cH}K) $ to be the single
complex associated with
$$
\aligned
K_{\bZ}\oplus (\Dec W)_{0}K_{\bQ}\oplus
F^{0}(\Dec W)_{0}K_{\bC}
&\to
K'_{\bQ}\oplus (\Dec W)_{0}K'_{\bC},
\\
K_{\bZ}\oplus K_{\bQ}\oplus F^{0}K_{\bC}
&\to
K'_{\bQ}\oplus K'_{\bC}
\endaligned
$$
respectively.
They can be defined also for a mixed Hodge complex
$ K $ in the sense of [5, 3.2] (where
$ \Dec W $ is replaced by
$ W $).
Using these, we can also define
$ H_{\cD}^{i}(X,A(k))' $,
$ H_{\cD}^{i}(X,A(k)) $.
Note that
$ \Gamma (D_{\cH}K) $ is canonically isomorphic to
$$
K_{\bZ}\oplus F^{0}K_{\bC} \to K'_{\bC},
\leqno(1.2.2)
$$
if there is a canonical morphism
$ \alpha' : K_{\bZ} \to K_{\bQ} $ such that
$ \alpha_{1} = \alpha_{2}\scirc\alpha' $.

\medskip\noindent
{\bf 1.3.~Lemma.} {\it
The canonical morphism
$ H_\cD^{i} (X,A(k))' \to H_\cD^{i} (X,A(k)) $ is an isomorphism if
$$
\text{
$ H^{i-1}(X,A) $ and
$ H^{i}(X,A) $ have weights
$ \le 2k $,
}
\leqno(1.3.1)
$$
and
$ H_\cD^{i} (X,A(k))'' \to H_\cD^{i} (X,A(k))' $ is an isomorphism if
$$
\text{
$ \Gr_{2k}^{W}H^{i-1}(X,\bQ) $ is isomorphic to a direct sum of
$ \bQ(-k) $.
}
\leqno(1.3.2)
$$
We have the corresponding assertion for Deligne homology where
$ H^{i-1}(X,A) $,
$ H^{i}(X,A) $ and
$ k $ are replaced respectively by
$ H_{i+1}^{\BM}(X,A) $,
$ H_{i}^{\BM}(X,A) $ and
$ - k $.
}

\medskip\noindent
{\it Proof.}
This is clear by (1.1.2).

\medskip\noindent
{\bf 1.4.~Remark.}
Condition (1.3.2) is satisfied for
$ H^{2p-2}(X \setminus Y,\bQ) $ with
$ i = 2p - 1 $ and
$ k = p $, if
$ X $ is smooth and
$ Y $ is a closed subvariety of codimension
$ \ge p - 1 $.
Indeed,
$ \Gr_{2p}^{W}H_{Y}^{2p-1}(X,\bQ) $ is a direct sum of
$ \bQ(-p) $.
A similar assertion holds also for
$ \Gr_{2-2m}^{W}H_{2m-1}^{\BM} (Y,\bQ) $ if
$ Y $ is of pure dimension
$ m $.

\medskip\noindent
{\bf 1.5.~Mixed Hodge Modules} (see [47]).
For a variety
$ X $ we denote by
$ \MHM(X) $ the abelian category of mixed
$ \bQ $-Hodge Modules on
$ X $,
and
$ D^{b}\MHM(X) $ its derived category consisting of bounded complexes of
mixed
$ \bQ $-Hodge Modules.
There is a natural functor
$ \rat : D^{b}\MHM(X) \to D_{c}^{b}(X,\bQ) $ assigning the underlying
$ \bQ $-complexes where
$ D_{c}^{b}(X,\bQ) $ denotes the full subcategory of
$ D_{c}^{b}(X^{\an},\bQ) $ consisting of
$ \bQ $-complexes whose cohomology sheaves are algebraically
constructible.
We denote by
$ H^{i} : D^{b}\MHM(X) \to \MHM(X) $ the usual cohomology functor.

For morphisms
$ f $ of algebraic varieties we have canonically defined functors
$ f_{*}, f_{!}, f^{*}, f^{!} $ between the derived categories of mixed
$ \bQ $-Hodge Modules.
They are compatible with the corresponding functors of
$ \bQ $-complexes via the functor
$ \rat $.
For a closed embedding
$ i : X \to Y $,
the direct image
$ i_{*} $ will be omitted sometimes in order to simplify the notation,
because
$$
i_{*} : D^{b}\MHM(X) \to D^{b}\MHM(Y)
\leqno(1.5.1)
$$
is fully faithful.

If
$ X = \Spec \, \bC $ we have naturally an equivalence of categories
$$
\MHM(\Spec \, \bC) = \MHS(\bQ)^{p}.
\leqno(1.5.2)
$$
Here the right-hand side is as in (1.1).
So
$ \MHM(\Spec \, \bC) $ will be identified with
$ \MHS(\bQ)^{p} $.

We denote by
$ \bQ(j) $ the mixed Hodge structure of type
$ (-j,-j) $ whose underlying
$ \bQ $-vector space is
$ (2\pi i)^{j}\bQ \subset \bC $, see [22].
For a variety
$ X $ with structure morphism
$ a_{X} : X \to \Spec \, \bC $,
we define
$$
\bQ_{X}^{H}(j) = a_{X}^{*}\bQ(j),\quad \bD_{X}^{H}(j) =
a_{X}^{!}\bQ(j),
\leqno(1.5.3)
$$
so that
$ \bD_{X}^{H}(j) $ is the dual of
$ \bQ_{X}^{H}(-j) $.
We will write
$ \bQ_{X}^{H} $ for
$ \bQ_{X}^{H}(0) $,
and similarly for
$ \bD_{X}^{H} $.
If
$ X $ is smooth of pure dimension
$ n $,
we have a canonical isomorphism
$$
\bD_{X}^{H} = \bQ_{X}^{H}(n)[2n].
\leqno(1.5.4)
$$

\medskip\noindent
{\bf 1.6.~Proposition.} {\it
With the notation of
$ (1.1) $ and
$ (1.5) $ we have canonical isomorphisms
$$
\aligned
H_{\cD}^{i}(X,\bQ(k))'' &= \Ext^{i}(\bQ,(a_{X})_{*}\bQ_{X}^{H}(k)),
\\
H_{i}^{\cD}(X,\bQ(k))'' &= \Ext^{-i}(\bQ,(a_{X})_{*}\bD_{X}^{H}(-k)).
\endaligned
$$
}

\medskip\noindent
{\it Proof.}
In the case
$ A = \bQ $, we have canonical isomorphisms
$ K = (a_{X})_{*}\bQ_{X}^{H} $,
$ K' = (a_{X})_{*}\bD_{X}^{H} $ by [50].

\medskip\noindent
{\bf 1.7.~Remark.}
Let
$ X $ be a reduced variety of pure dimension
$ n $,
and
$ X_{i} $ be the irreducible components of
$ X $.
Let
$ \Rat(X)^{*} = \prod \Rat(X_{i})^{*} $ with
$ \Rat(X_{i}) $ the rational function field of
$ X_{i} $.
Then by [27, 2.12], [32, 3.1] we have a canonical isomorphism
$$
H_{2n-1}^{\cD}(X,\bZ(n-1)) = \{g \in \Rat(X)^{*} : \div g = 0\},
\leqno(1.7.1)
$$
where
$ \div g = \sum \div g_{i} $ if
$ g = (g_{i}) $ with
$ g_{i} \in \Rat(X_{i})^{*} $.
(Here (1.3.1--2) are satisfied.)

If
$ X $ is smooth, this is due to [27, 2.12].
In this case, the left-hand side of (1.7.1) is isomorphic to
$ \Ext^{1}(\bZ_{X},\bZ_{X}(1)) $ (where
$ \Ext^{1} $ is taken in the category of admissible variation of mixed
Hodge structures), and the assertion is related with the theory of
$ 1 $-motives [22], and is more or less well-known.
Indeed, if
$ X $ is a point, the assertion is verified by calculating the period of
the mixed Hodge structure on
$ H^{1}(\bA^{1}\setminus \{0\}, \{1\}\cup \{x\}) $ for
$ x \in \bC \setminus \{0,1\} $, i.e., by using the integral of
$ dt/t $ on the relative cycle connecting
$ \{1\} $ and
$ \{x\} $,
where
$ t $ is the coordinate of
$ \bA^{1} $.
The general case is reduced to the smooth case using a long exact
sequence, see [32, 3.1].

\bigskip\bigskip
\centerline{{\bf 2. Higher Chow groups and cycle maps}}

\bigskip
\noindent
{\bf 2.1.~Higher Chow groups} ([10]).
Let
$ \Delta^{n} = \Spec (\bC[t_{0}, \dots, t_{n}]/(\sum t_{i} - 1)) $.
For a subset
$ I $ of
$ \{0, \dots, n\} $,
let
$ \Delta_{I}^{n} = \{t_{i} = 0 \, (i \in I)\} \subset \Delta^{n} $.
It is naturally isomorphic to
$ \Delta^{m} $ with
$ m = n - |I| $ (fixing the order of the coordinates), and is called a
face of
$ \Delta^{n} $.
For
$ 0 \le i \le n $,
we have inclusions
$ \iota_{i} : \Delta^{n-1} \to \Delta^{n} $ such that
its image is
$ \Delta_{\{i\}}^{n} $.

Let
$ X $ be an equidimensional variety.
Then
$ X \times \Delta_{I}^{n} $ is also called a face of
$ X \times \Delta^{n} $.
Following Bloch, we define
$ \cZ^{p}(X,n) $ to be the free abelian group with generators the
irreducible closed subvarieties of
$ X \times \Delta^{n} $ of codimension
$ p $,
intersecting all the faces of
$ X \times \Delta^{n} $ properly.
We have face maps
$$
\partial_{i} : \cZ^{p}(X,n) \to \cZ^{p}(X,n-1),
$$
induced by
$ \iota_{i} $.
Let
$ \partial = \sum (-1)^{i}\partial_{i} $.
Then
$ \partial^{2} = 0 $,
and
$ \CH^{p}(X,n) $ is defined to be
$ \Ker \, \partial /\Im \, \partial $ which is a subquotient of
$ \cZ^{p}(X,n) $.
By [10] it is isomorphic to
$$
\frac{{\bigcap}_{0\le i\le n}\Ker({\partial}_{i}:\cZ^{p}(X,n)\to
\cZ^{p}(X,n-1))}
{{\partial}_{n+1}({\bigcap}_{0\le i\le n}
\Ker({\partial}_{i}:\cZ^{p}(X,n+1)\to \cZ^{p}(X,n)))}
\leqno(2.1.1)
$$

Indeed, let
$ \cZ^{p}(X,\cssbull)' $ be the subcomplex of
$ \cZ^{p}(X,\cssbull) $ defined by
$$
\cZ^{p}(X,n)' = {\bigcap}_{0\le i< n}
\Ker({\partial}_{i}:\cZ^{p}(X,n)\to \cZ^{p}(X,n-1)).
$$
Then the inclusion induces a quasi-isomorphism
$$
\cZ^{p}(X,\cssbull)' \to \cZ^{p}(X,\cssbull).
\leqno(2.1.2)
$$
(For this, we can consider first the subcomplex defined by
$ \Ker \, \partial_{0} $,
using a homotopy given by the zeroth degeneracy, and then proceed
inductively.)

\medskip\noindent
{\bf 2.2.~Remarks.}
(i)
In this paper we are mainly interested in
$ \CH^{p}(X,n) $ for
$ n = 0, 1 $.
If
$ n = 0 $, it is the usual Chow group.
If
$ n = 1 $, any higher cycle
$ \zeta \in \CH^{p}(X,1) $ can be represented by
$ \msum_{j} (Z_{j},g_{j}) $ where
$ Z_{j} $ are irreducible (and reduced) subvarieties of
$ X $ with pure codimension
$ p - 1 $ and
$ g_{j} $ are rational functions on
$ Z_{j} $ such that
$ \msum_{j} \div g_{j} = 0 $.
Indeed, such elements modulo the relation given by the tame symbols form
an abelian group
$ H^{p-1}(X,\cK_{p}) $ using the Gersten resolution, where
$ \cK_{p} $ is the Zariski-sheafification of the Quillen
$ K $-group.
It is well known (see e.g. [40]) that there is a natural isomorphism
$$
H^{p-1}(X,\cK_{p}) = \CH^{p}(X,1).
\leqno(2.2.1)
$$
For each
$ \msum_{j} (Z_{j},g_{j}) $, the corresponding higher cycle is defined by
taking the closure of the graph of
$ g_{j} $ in
$ X \times \bP^{1} $, and then restricting it to the complement of
$ X \times \{1\} $.
Here we use an automorphism of
$ \bP^{1} $ sending
$ 0, 1, \infty $ to
$ 0, \infty, 1 $ respectively
(or rather take another affine coordinate of
$ \bP^{1} $).

\medskip
(ii)
If we assume
$ g_{j} = \text{\rm const} $ in the above Remark, we get a natural
morphism
$$
\CH^{p-1}(X)\otimes \bC^{*} \to \CH^{p}(X,1).
\leqno(2.2.2)
$$
Its image is denoted by
$ \CH_{\dec}^{p-1}(X,1) $, and is called the subgroup of
{\it decomposable} higher cycles, see [19], [40], etc.
We define the group of indecomposable higher cycles by
$$
\CH_{\ind}^{p-1}(X,1)_{\bQ} =
\CH^{p-1}(X,1)_{\bQ}/\CH_{\dec}^{p-1}(X,1)_{\bQ}.
\leqno(2.2.3)
$$

\medskip\noindent
{\bf 2.3.~Functoriality.} Let
$ f : X \to Y $ be a proper morphism of varieties, and put
$ r = \dim X - \dim Y $.
Then we have the pushforward functor
$$
f_{*} : \CH^{p}(X,n) \to \CH^{p-r}(Y,n).
$$
In fact, for a face map
$ \iota : \Delta^{m} \to \Delta^{n} $,
Bloch showed the commutative diagram
$$
\CD
\cZ^{p}(X,n) @>{{\iota}^{*}}>> \cZ^{p}(X,m)
\\
@V{f_{*}}VV @V{f_{*}}VV
\\
\cZ^{p-r}(X,n) @>{{\iota}^{*}}>> \cZ^{p-r}(X,m)
\endCD
$$

As for the pull-back, we have
$ f^{*} : \CH^{p}(Y,n) \to \CH^{p}(X,n) $ if
$ f $ is flat.
In the case
$ X, Y $ are quasi-projective and smooth, we have
$ f^{*} : \CH^{p}(Y,n)_{\bQ} \to \CH^{p}(X,n)_{\bQ} $ by
[38].
Here we have a quasi-isomorphic subcomplex
$ \cZ_{f}^{p}(Y,\cssbull)_{\bQ} $ of
$ \cZ^{p}(Y,\cssbull)_{\bQ} $ on which the pull-back
$ f^{*} $ is naturally defined.

\medskip\noindent
{\bf 2.4.~Cycle map.} Let
$ X $ be an equidimensional variety.
By [4], [11], [23], etc., we have a cycle map
$$
cl : \CH^{p}(X,n) \to H_{2d+n}^{\cD}(X,\bQ(d))'',
\leqno(2.4.1)
$$
where
$ d = \dim X - p $.
The target becomes
$ H_{\cD}^{2p-n}(X,\bQ(p))'' $ by (1.1.1) if
$ X $ is smooth.
Using mixed Hodge Modules [48], the cycle map (2.4.1) is defined
as follows.

Let
$ S^{n-1} = \bigcup \Delta_{\{i\}}^{n} \subset \Delta^{n}, U =
\Delta^{n} \setminus S^{n-1} $ with the inclusion morphisms
$ i : S^{n-1} \to \Delta^{n}, j : U \to \Delta^{n} $.
Then
$$
(a_{\Delta^{n}})_{*}j_{!}\bQ_{U}^{H} = \bQ_{pt}^{H}[-n],
\leqno(2.4.2)
$$
where
$ a_{\Delta^{n}} : \Delta^{n} \to pt := \Spec \, \bC $ is the
structure morphism.
Let
$ \zeta = \msum_{k} n_{k}[Z_{k}] \in \bigcap_{0\le i\le n} \Ker \,
\partial_{i}\subset \cZ^{p}(X,n) $ (see (2.1)), where
$ Z_{k} $ are irreducible closed subvarieties of
$ X \times \Delta^{n} $.
Let
$ d' = \dim Z_{k} = d + n $.
Put
$ Z = \bigcup_{k} Z_{k} $.
Then the coefficients
$ n_{k} $ of
$ Z_{k} $ induce a morphism
$$
\bQ_{Z}^{H} \to \moplus_{k}\, \IC_{Z_{k}}\bQ^{H}[-d']
\to \bD_{Z}^{H}(-d')[-2d'] \to
\bD_{X\times \Delta^{n}}^{H}(-d')[-2d'],
\leqno(2.4.3)
$$
where
$ \IC_{Z_{k}}\bQ^{H} $ denotes the mixed Hodge Module whose underlying
perverse sheaf is the intersection complex
$ \IC_{Z_{k}}\bQ $ [7].
Let
$ \pi : X \times \Delta^{n} \to X $ be the first projection, and
let
$ j $ denote also
$ \id \times j : X \times U \to X \times \Delta^{n} $ (and the same for
$ i) $.
Then
$$
\pi_{*}j_{!}\bD_{X\times U}^{H} = \bD_{X}^{H}(n)[n]
$$
by (1.5.4) and (2.4.2).
So it is enough to show that (2.4.3) is uniquely lifted to
$$
\bQ_{Z}^{H} \to j_{!}\bD_{X\times U}^{H}(-d')[-2d'],
$$
i.e., the composition of (2.4.3) with
$$
\bD_{X\times \Delta^{n}}^{H}(-d')[-2d'] \to
i_{*}i^{*}\bD_{X\times \Delta^{n}}^{H}(-d')[-2d']
$$
is zero and
$ \Hom(\bQ_{Z}^{H}, i_{*}i^{*}\bD_{X\times \Delta^{n}}^{H}(-d')[-2d'-1]) =
0 $.
But they follow from the condition on proper intersection together with
$ \zeta \in \bigcap_{0\le i\le n} \Ker \, \partial_{i} $.
For the well-definedness of the cycle map, it is enough to show its
invariance under a deformation of cycle parametrized by
$ \bA^{1} $ (using a blow-up of
$ \Delta^{n} $).

\medskip\noindent
{\bf 2.5~Remarks.} (i)
The cycle map for
$ n = 0 $ is defined with integral coefficients by the composition of
$$
\bZ \to H_{2d}^\cD(Z,\bZ(d))
\to H_{2d}^\cD(X,\bZ(d))
\leqno(2.5.1)
$$
for
$ \zeta = \msum_{k} n_{k}[Z_{k}] \in \CH_{d}(X) $, where
$ Z = \cup_{k} Z_{k} $, see [4], [29], etc.
If
$ X $ is smooth proper, this coincides with Deligne's cycle map,
which is defined by the composition of
$$
\bZ \to H_{Z}^{2p}(X,\bZ(p))
\to K(p)[2p]
$$
where
$ K $ is as in (1.1).
It induces Griffiths' Abel-Jacobi map
$$
\CH_{\hom}^{p}(X) \to J^{p}(X)
:= \Ext_{\MHS}^{1}(\bZ, H^{2p-1}(X,\bZ)(p)),
\leqno(2.5.2)
$$
see [25] (and also [47, (4.5.20)], etc.)
Here
$ J^{p}(X) $ is Griffiths' intermediate Jacobian by [16].
(This can be defined even if
$ X $ is not proper.)
If
$ p = \dim X $, then (2.5.2) is the Albanese map.

\medskip
(ii)
If
$ X $ is purely
$ m $-dimensional, the cycle map induces isomorphisms for
$ n > 0 $:
$$
\aligned
&\CH^{1}(X,n) \simto H_{2m-2+n}^{\cD}(X,\bZ(m-1))^{\prime\prime}
\\
&\qquad = H_{2m-2+n}^{\cD}(X,\bZ(m-1))^{\prime}
= H_{2m-2+n}^{\cD}(X,\bZ(m-1))
\endaligned
\leqno(2.5.3)
$$
and these groups vanish if
$ n > 1 $.
Indeed, let
$ U $ be a smooth dense open subvariety of
$ X $, and
$ Z $ its complement.
Then
$ \CH^{0}(Z,n) = 0 $ for
$ n > 0 $ by (2.1.1), and
$ \CH^{1}(U,n) = 0 $ for
$ n > 1 $ by [10, 6.1].
Thus
$ \CH^{1}(X,n) = 0 $ for
$ n > 1 $ by the localization sequence [12].
On the other hand, the last two isomorphisms of (2.5.3) follow from
(1.3), and
$ H_{2m-2+n}^{\cD}(X,\bZ(m-1))'' $ vanishes for
$ n > 1 $ by (1.1.3).
So the case
$ n > 1 $ is clear.
The case
$ n = 1 $ is shown by [27, 2.12] in the smooth case and [32, 3.1]
(combined with the above (2.2.1)) in general.
We will generalize (2.5.3) to the case
$ n = 0 $ in (3.1).

\medskip\noindent
{\bf 2.6.~Compatibility.} The cycle map (2.4.1) is compatible with
$ f_{*} $ for a proper morphism
$ f $, and also with
$ f^{*} $ for a morphism of smooth quasi-projective varieties
$ f : X \to Y $ in the case of rational coefficients.
Indeed, this is reduced to the case of the usual Chow groups by
(2.3) and the construction of the cycle map (2.4), and follows from [49].

\medskip\noindent
{\bf 2.7.~Proposition.} {\it
Let
$ X $ be a quasi-projective variety, and
$ Y $ a closed subvariety.
Assume
$ X, Y $ are equidimensional.
Let
$ r = \codim_{X} Y $,
and
$ d = \dim X - p $.
Then the cycle map induces a morphism of long exact sequences
$$
\CD
\CH^{p}(X\setminus Y,n+1) @>>> \CH^{p-r}(Y,n) @>>>
\CH^{p}(X,n)
\\
@VVV @VVV @VVV
\\
H_{2d+n+1}^\cD(X\setminus Y,\bQ(d))'' @>>>
H_{2d+n}^\cD(Y,\bQ(d))'' @>>> H_{2d+n}^\cD(X,\bQ(d))''
\endCD
$$
where the first exact sequence is the localization sequence
$ [12] $ {\rm (}choosing the sign appropriately{\rm ),}
and the second comes from
$ (1.1.4) $.
}

\medskip\noindent
{\it Proof.}
The assertion is clear except for the commutativity of the left part of
the diagram.
Let
$ U' = \Delta^{n+1} \setminus \cup_{0\le i\le n}
\Delta_{\{i\}}^{n+1} $,
and identify
$ \Delta_{\{n+1\}}^{n+1} \cap U' $ with
$ U^{n} := \Delta^{n} \setminus S^{n-1} $.
Let
$ j : X \setminus Y \to X, j' : U' \to \Delta^{n+1},
j^{n} : U^{n} \to \Delta^{n} $ denote the inclusion morphisms so
that we have distinguished triangles
$$
\gathered
\to \bD_{Y}^{H} \to \bD_{X}^{H}
\to j_{*}\bD_{X\setminus Y}^{H}
\to,
\\
\to j_{!}^{n+1}\bD_{{U}^{n+1}}^{H}
\to j'_{!}\bD_{U'}^{H} \to
j_{!}^{n}\bD_{{U}^{n}}^{H}(1)[2] \to,
\endgathered
$$
where the direct images by closed embeddings are omitted to simplify
the notation, see (1.5.1).

Let
$ \zeta \in \CH^{p}(X\setminus Y,n+1) $.
By (2.1.2) and [12], it is represented by
$ \zeta \in \cZ^{p}(X\setminus Y,n+1)' $ which is extended to
$ \zeta' \in \cZ^{p}(X,n+1)' $ so that its restriction to
$ X\times \Delta^{n} $ is
$ \oze \in \cZ^{p}(Y,n)' $. Then the image of
$ \zeta $ by the morphism of the localization sequence is
$ \oze $.
By definition
$ \zeta' $ gives
$$
\xi \in \Hom(\bQ_{Z}^{H},
j'_{!}\bD_{X\times U'}^{H} (-d'-1)[-2d'-2])
$$
such that its restriction to
$ (X \setminus Y) \times \Delta^{n} $ vanishes.
(Here
$ d' = \dim X -p + n $,
and
$ j' $ denotes also
$ id \times j' $.)
So it induces
$$
\aligned
\xi'
&\in \Hom(\bQ_{Z}^{H},
j_{!}^{n}\bD_{Y\times {U}^{n}}^{H} (-d')[-2d']),
\\
\xi''
&\in \Hom(\bQ_{Z}^{H}, j_{*}j_{!}^{n+1}
\bD_{(X\setminus Y)\times {U}^{n+1}}^{H} (-d'-1)[-2d'-2]),
\endaligned
$$
using the external product of the above two distinguished triangles.
We see that
$ \xi' $ coincides with the image of
$ \oze $ by the cycle map, and the second distinguished triangle
induces an isomorphism
$ (a_{\Delta^{n}})_{*}j_{!}^{n}\bD_{{U}^{n}}^{H} (1)[2]
\to (a_{\Delta^{n+1}})_{*}j_{!}^{n+1}\bD_{{U}^{n+1}}^{H} $.
So the assertion is reduced to the next lemma.
(Here we represent the middle terms of the distinguished triangles by the
mapping cone of the morphism of the other terms so that we get short exact
sequences as below.)

\medskip\noindent
{\bf 2.8.~Lemma.} {\it
Let
$ \{K^{i,j,\ssbull}\} $ be a square diagram of short exact sequences of
complexes of an abelian category, i.e.
$ K^{i,j,k} = 0 $ for
$ |i| > 1 $ or
$ |j| > 1 $, and
$ K^{i-1,j,k} \to K^{i,j,k} \to K^{i+1,j,k} $ is exact
{\rm (}and similarly for the index
$ j )$.
Assume
$ H^{k-1}(K^{1,1,\ssbull}) = 0 $.
Let
$ \xi \in H^{k}(K^{0,0,\ssbull}) $ such that its image in
$ H^{k}(K^{1,1,\ssbull}) $ vanishes.
Let
$ \xi' \in H^{k}(K^{-1,1,\ssbull}) $,
$ \xi'' \in H^{k}(K^{1,-1,\ssbull}) $ such that the images of
$ \xi $,
$ \xi' $ in
$ H^{k}(K^{0,1,\ssbull}) $ coincide and the images of
$ \xi $,
$ \xi'' $ in
$ H^{k}(K^{1,0,\ssbull}) $ coincide.
Then the images of
$ \xi' $,
$ \xi'' $ in
$ H^{k+1}(K^{-1,-1,\ssbull}) $ coincide up to sign.
}

\medskip
(The proof is straightforward.
For a similar assertion, where the triangle is slightly
shifted, a proof is given in [34], p. 268.)

\bigskip\bigskip
\centerline{{\bf 3. Proof of main theorems and related assertions}}

\medskip\noindent
Since Theorem (0.2) follows from (0.3) and (2.7),
we first show Theorem (0.3).

\medskip\noindent
{\bf 3.1.~Proof of Theorem (0.3).}
It is enough to show the assertion for
$ H_{2m-2}^\cD (Y, \bZ(m-1))'' $ by (1.3).
We apply (2.7) to
$ Y $ and a divisor
$ Z $ on $ Y $ containing
$ \Sing Y $.
Let
$ U = Y \setminus Z $.
The assertion follows from [49, I, (3.4)] if
$ Y $ is smooth (i.e. if
$ U = Y $).
Note that we have the surjectivity of the cycle map
$ \CH^{1}(U) \to H_{2m-2}^\cD(U, \bZ(m-1))'' $ in
loc. cit, because
$ H^{1}(U,\bZ) $ is torsion-free.
By [27, 2.12], we have a similar isomorphism
$$
\CH^{1}(U,1) \to H_{2m-1}^{\cD}(U,\bZ(m-1))''.
\leqno(3.1.1)
$$
So the general case is reduced to the smooth case
by using the cycle map of the localization sequence
$$
\CH^{1}(U,1) \to \CH^{0}(Z) \to
\CH^{1}(Y) \to \CH^{1}(U) \to 0
$$
to the corresponding exact sequence of Deligne homology.
Indeed, the cycle maps are compatible with the first morphism
$ \CH^{1}(U,1) \to \CH^{0}(Z) $ (which is given by the divisor map) by
using (1.7.1) and [10, 6.1], see [32, 3.1].

\medskip\noindent
{\bf 3.2.~Remark.}
In general, the isomorphism
$ \CH^{1}(Y) = H_{2m-2}^\cD(Y,\bZ(m-1)) $ does not hold
(even for a smooth
$ Y $), see [49, I, (3.5)].
(In Remark (i) of loc. cit. the assumption of the
second statement should be replaced by the condition that
$ H^{2}(X,\bQ) \cap F^{1}H^{2}(X,\bC) $ is not contained in
$ W_{2}H^{2}(X,\bQ) $.)

\medskip\noindent
{\bf 3.3.~Theorem.} {\it
Let
$ X $ be a connected smooth projective variety.
For an open subvariety
$ U $ of
$ X,$ let
$ j_{U,X} : U \to X $ be the inclusion morphism, and
$ j_{U,X}^{*} $ the pull-back of Deligne cohomology.
Then, for any integer
$ p $,
the following three conditions are equivalent to each other:

\smallskip\noindent
$ (a) $
Griffiths' Abel-Jacobi map
$ \CH_{\hom}^{p}(X)_{\bQ} \to J^{p}(X)_{\bQ} $ is injective.

\smallskip\noindent
$ (b) $
$ \underset{Y}\to\varinjlim\, H_\cD^{2p-1}(X \setminus Y,\bQ(p))'/
\Im \,j_{X\setminus Y,X}^{*}=0 $,
where
$ Y $ runs over the closed subvarieties of
$ X $ with pure codimension
$ p - 1 $.

\smallskip\noindent
$ (c) $
The cycle map
$ \CH^{p}(X \setminus Y,1)_{\bQ} \to H_\cD^{2p-1}(X \setminus
Y,\bQ(p))'/\Im \, j_{X\setminus Y,X}^{*} $ is surjective for any
{\rm (}sufficiently large{\rm )} closed subvarieties
$ Y $ of
$ X $ with pure codimension
$ p - 1 $.
}

\medskip\noindent
{\it Proof.}
The equivalence of (b) and (c) follows from (2.7)
together with (2.5.3) for
$ n = 1 $, which we apply to closed subvarieties of pure codimension
$ p - 1 $ in
$ X $.
The equivalence of (a) and (c) follows from (2.7) and (0.3).

\medskip\noindent
{\bf 3.4.~Remarks.}
(i) We have
$ H_\cD^{2p-1}(X \setminus Y,\bQ(p))' = H_\cD^{2p-1}(X \setminus Y,
\bQ(p))'' $ if
$ Y $ has codimension
$ \ge p - 1 $,
and
$ H_\cD^{2p-1}(X \setminus Y,\bQ(p)) = H_\cD^{2p-1}(X \setminus Y,
\bQ(p))' $ if furthermore
$ p = \dim X $.
(In this case, condition (a) is related with [44].)
The pull-back
$ j_{X \setminus Y,X}^{*} $ vanishes if
$ H^{2p-2}(X,\bQ)(p-1) $ is generated by algebraic cycle classes and
$ Y $ is sufficiently large (using (1.1.4)).
I am informed that the three conditions (a), (b), (c)
are further equivalent to the condition:

\smallskip\noindent
$ (d) $
The cycle map
$ \CH^{p}(X \setminus Y,1)_{\bQ} \to
\Hom_{\MHS}(\bQ,H^{2p-1}(X \setminus Y,\bQ(p)) $ is surjective
for any (sufficiently large) closed subvarieties
$ Y $ of
$ X $ with pure codimension
$ p $.

\smallskip\noindent
This was studied by Jannsen ([33], 9.10), and is equivalent to
condition (c) with last
$ p - 1 $ replaced by
$ p $, see 9.8 in loc. cit.

\medskip
(ii) The exact sequence of L. Barbieri-Viale and V. Srinivas [2]
for cycles of codimension
$ 2 $ in the introduction follows also from the local-to-global
spectral sequence in the theory of Bloch and Ogus [14] applied
to the absolute Hodge cohomology (using (0.3) in the smooth case).
Here the flasqueness of
$ \cH_{\cD}^{3}(\bZ(2))'' $ is clear, because the inductive
limit of
$ H_{\cD}^{2}(U,\bZ(1))'' $ vanishes where
$ U $ runs over the nonempty open subvarieties of an irreducible
subvariety of codimension
$ 1 $ in
$ X $.

\medskip\noindent
{\bf 3.5.~Relation with the normalization.}
We can express the subgroup
$ \CH_{\hom}^{1}(Y) $ of
$ \CH^{1}(Y) $ consisting of Borel-Moore homologically equivalent to zero
cycles by using the normalization of
$ Y $.
This may be useful for explicit calculation.

Let
$ Y $ be a connected variety of pure dimension
$ m $ with
$ Y_{k}\, (1 \le k \le r) $ the irreducible components of
$ Y $.
Let
$ \tY $ be the disjoint union of the normalizations
$ \tY_{k} $ of
$ Y_{k} $ with
$ \pi : \tY \to Y $ the natural morphism.
Let
$ D = \{y \in Y : |\pi^{-1}(y)| > 1\} $.
We assume
$ \tY $ is smooth,
$ D $ is a smooth closed subvariety of
$ Y $ with pure codimension one, and
$ \pi_{*}\bZ_{\tY}|_{D} $ is a local system.
(We may assume these because
$ \CH^{1}(Y) $ and
$ H_{2m-2}^\cD(Y,\bZ(m-1)) $ do not change by deleting a
closed subvariety of codimension
$ > 1 $.)

Let
$ \tD = \pi^{-1}(D) $,
and
$ \tD_{i} \,( i \in I), D_{j} \,( j \in J) $ be connected components of
$ \tD, D $.
Put
$$
I_{j} = \{i \in I : \tD_{i} \subset \pi^{-1}(D_{j})\},\quad
I(k) = \{i \in I : \tD_{i} \subset \tY_{k}\}.
$$
We define
$ \cE_{j} = \Ker(\Tr : \moplus_{i\in I_{j}}
\pi_{*}\bZ_{\tD_{i}} \to \bZ_{D_{j}}) $,
$ E_{j} = H^{0}(D_{j}, \cE_{j}) $, and
$ E = \moplus_{j} E_{j} $.
Let
$ d_{i} $ be the degree of
$ \tD_{i} $ over
$ \pi (\tD_{i}) $.
Then
$ E_{j} $ is naturally identified with
$$
\{a_{i} \in \bZ \,(i \in I_{j}) : \msum_{i\in I_{j}} d_{i}a_{i} = 0\}.
$$
Let
$ E'= \moplus_{j\in J} H^{1}(D_{j},\cE_{j}) $.
(This may have torsion which is related to the cokernel of
$ \CH^{1}(\tY) \to \CH^{1}(Y) $.)
Then we have an exact sequence
$$
\aligned
0 \to H_{2m-1}^{\BM}(\tY,\bZ)
&\to H_{2m-1}^{\BM}(Y,\bZ)
\to E(m-1)
\\
\overset\gamma\to\to
H_{2m-2}^{\BM}(\tY,\bZ)
&\to H_{2m-2}^{\BM}(Y,\bZ)
\to E'(m-1),
\endaligned
\leqno(3.5.1)
$$
where
$ \gamma $ is defined by
$ (a_{i}) \to \msum_{i} a_{i} cl([\tD_{i}]) $.
Here
$ cl([\tD_{i}]) $ denotes the cycle class.

We define
$ E^{0} = \Ker \,\gamma \subset E $ so that we get
$$
0 \to H^{1}(\tY,\bZ) (1)
\to H_{2m-1}^{\BM}(Y,\bZ)(1-m)
\to E^{0} \to 0,
\leqno(3.5.2)
$$
The associated extension class is denoted by
$ e \in \Ext_{\MHS}^{1}(E^{0}, H^{1}(\tY,\bZ)(1)) $.
The cycle map induces an isomorphism of exact sequences
$$
\CD
E^{0} @>>> \CH_{\hom}^{1}(\tY) @>>> \CH_{\hom}^{1}(Y) @>>> 0
\\
@| @| @|
\\
E^{0} @>>> J^{1}(\tY) @>>> J^{1}(Y)^{\BM} @>>> 0,
\endCD
\leqno(3.5.3)
$$
where
$ J^{1}(Y)^{\BM} := \Ext_{\MHS}^{1}(\bZ,
H_{2m-1}^{\BM}(Y,\bZ)(1-m)) $ and
$ J^{1}(\tY) $ is as in (2.5.2)
(and is a quotient of the Jacobian of a smooth compactification of
$ \tY $).

Indeed, let
$ \CH^{1}(Y)' = \Im(\CH^{1}(\tY)
\to \CH^{1}(Y)) $,
$ \CH_{\hom}^{1}(Y)^{\prime} = \CH^{1}(Y)' \cap \CH_{\hom}^{1}(Y) $.
Then, for the exactness of the first row, it is sufficient to show
$$
\CH_{\hom}^{1}(Y)^{\prime} = \CH_{\hom}^{1}(Y).
\leqno(3.5.4)
$$
This is reduced to the case where the cycle is supported on
$ \Sing Y $, and follows from the localization sequence for
Borel-Moore homology.
The second row is induced by (3.5.2), and we can show that for
$ u \in \Hom_{\MHS}(\bZ,E^{0}) $,
$$
e\scirc u \in
\Ext_{\MHS}^{1}(\bZ, H^{1}(\tY,\bZ)(1))
\leqno(3.5.5)
$$
coincides with the image of
$ \msum_{i} a_{i} [\tD_{i}] $ by the Abel-Jacobi map,
where
$ (a_{i}) = u(1) \in E^{0} $.
This is verified by using a natural morphism of (3.5.2) to
$$
0 \to H^{1}(\tY,\bZ)(1) \to
H^{1}(\tY \setminus \tD,\bZ)(1) \to
H^{0}(\tD,\bZ).
$$

\medskip\noindent
{\bf 3.6.~Griffiths group.} Let
$ \CH_{\alg}^{p}(X) $ denote the subgroup of cycles algebraically
equivalent to zero, and
$ \CH_{\AJ}^{p}(X) $ denote the kernel of Griffiths' Abel-Jacobi
map [30].
Let
$ J^{p}(X)^{\alg} $ be the image of the Abel-Jacobi map, and
$ J^{p}(X)^{\ab} $ the image of
$ \CH_{\alg}^{p}(X) $ which is an abelian subvariety of
$ J^{p}(X) $.
We call
$ J^{p}(X)^{\ab} $ and
$ J^{p}(X)^{\alg}/J^{p}(X)^{\ab} $ respectively the abelian and
discrete part of the image of the Abel-Jacobi map.
Let
$ \Griff^{p}(X) $ denote the Griffiths group
$ \CH_{\hom}^{p}(X)/\CH_{\alg}^{p}(X) $,
and
$ \Griff_{\AJ}^{p}(X) $ the image of
$ \CH_{\AJ}^{p}(X) $ in
$ \Griff^{p}(X) $.
Then we have a commutative diagram of short exact sequences
(where the
$ 0 $ are omitted):
$$
\CD
\CH_{\alg}^{p}(X) \cap \CH_{\AJ}^{p}(X) @>>>
\CH_{\alg}^{p}(X) @>>> J^{p}(X)^{\ab}
\\
@VVV @VVV @VVV
\\
\CH_{\AJ}^{p}(X) @>>> \CH_{\hom}^{p}(X) @>>> J^{p}(X)^{\alg}
\\
@VVV @VVV @VVV
\\
\Griff_{\AJ}^{p}(X) @>>> \Griff^{p}(X) @>>>
J^{p}(X)^{\alg} /J^{p}(X)^{\ab}
\endCD
$$

It is known that
$ \CH_{\alg}^{p}(X) $ is divisible [13].
By [39], [45] we have

\medskip\noindent
(3.6.1)\,\,\,\,
$ \CH_{\AJ}^{p}(X) $ is torsion-free if
$ p = 2 $ or
$ \dim X $.

\medskip\noindent
For
$ p = 2 $, this is proved in [39] by using Bloch's cycle map [8].
(See also [53]).

\medskip\noindent
{\bf 3.7~Remark.}
For an open subvariety
$ U $ of
$ X $ and a positive integer
$ m $, consider the short exact sequence
$$
0 \to H_{\cD}^{2}(U,\bZ(2))/m \to H^{2}(U,\bZ(2)/m) \to
{}_{m}(H_{\cD}^{3}(U,\bZ(2))) \to 0,
\leqno(3.7.1)
$$
where
$ M/m = M/mM $ and
$ {}_{m}M = \Ker(m : M \to M )$ for an abelian group
$ M $.
The first morphism is surjective after taking the inductive
limit over
$ U $ by [39], and the last term is related to the kernel of the
Abel-Jacobi map, see (0.2).
So this can be used also for the proof of (3.6.1).

\medskip\noindent
{\bf 3.8.~Proposition.}
{\it If
$ p = 2 $,
$ \CH_{\alg}^{2}(X) \cap \CH_{\AJ}^{2}(X) $ is divisible and
torsion-free {\rm (}i.e. a
$ \bQ $-vector space{\rm ).}
}

\medskip\noindent
{\it Proof.}
Since the torsion-freeness follows from (3.6.1), it is enough
to show the divisibility.
Let
$ \zeta \in \CH_{\alg}^{2}(X) \cap \CH_{\AJ}^{2}(X) $.
This comes from
$ \zeta' \in \Pic^{0}(Y) $ where
$ Y $ is a resolution of singularities of a closed subvariety of
pure codimension
$ 1 $ in
$ X $.
Let
$ P = \Ker(\Pic^{0}(Y) \to J^{2}(X)) $.
Then
$ P $ is an extension of a finite group
$ \Gamma $ by an abelian variety
$ P^{0} $.
Since
$ P^{0} $ is divisible, there is an exact sequence
$$
0 \to {P}_{\tor}^{0} \to P_{\tor} \to \Gamma \to 0,
$$
and it splits.
Let
$ P' = \Ker(\Pic^{0}(Y) \to \CH_{\alg}^{2}(X)) $.
Since
$ P/P' $ is torsion-free by (3.6.1), we get an isomorphism
$ P'_{\tor} \to P_{\tor} $,
and
$ P^{0} \to P/P' $ is surjective, i.e.
$$
P^{0}/(P^{0} \cap P') = P/P'.
$$
So the assertion follows, since the left-hand side is divisible.

\medskip\noindent
{\bf 3.9.~Corollary.}
{\it If
$ p = 2 $,
$ \Griff_{\AJ}^{2}(X) $ is torsion-free.
}

\medskip\noindent
{\it Proof.}
This is clear by (3.8) and (3.6.1),
using the left column of the commutative diagram in (3.6).

\medskip\noindent
{\bf 3.10.~Proposition.}
{\it Let
$ X $ be an irreducible smooth proper complex algebraic
variety with a surjective morphism
$ f : X \to S $ whose general fibers
$ X_{s} $ are connected and
have no global
$ 2 $-forms.
If general fibers have dimension
$ > 2 $, we assume that the monodromy invariant
part of
$ H^{3}(X_{s},\bQ) $ vanishes by restricting it to a
sufficiently small nonempty open subvariety of
$ X_{s} $ for a general
$ s \in S $.
If
$ \dim S > 1 $, we assume further that
$ H^{1}(X_{s},\bQ) = 0 $ for a general
$ s \in S $, and the Abel-Jacobi map is surjective for cycles
of codimension
$ 2 $ on
$ S $.
Then the Abel-Jacobi map
$ \CH_{\hom}^{2}(X) \to J^{2}(X) $ is surjective
and the discrete part
$ J^{2}(X)^{\alg} /J^{2}(X)^{\ab} $ vanishes.
}

\medskip\noindent
{\it Proof.}
Since
$ J^{p}(X)^{\alg}/J^{p}(X)^{\ab} $ is discrete,
$ J^{p}(X) = J^{p}(X)^{\alg} $ if and only if
$ J^{p}(X) = J^{p}(X)^{\ab} $.
So the assertion is equivalent to the vanishing of
$ \Gr_{3}^{W}H^{3}(U,\bQ) $ for a sufficiently small non-empty
affine open subvariety
$ U $.
Indeed, let
$ Y $ be a resolution of singularities of the divisor
$ X \setminus U $.
Then, using the localization sequence, these two conditions
are both equivalent to the surjectivity of
$ H^{1}(Y,\bQ) \to H^{3}(X,\bQ)(1) $, or equivalently, of
$ J^{1}(Y) \to J^{2}(X) $.

To show the vanishing of
$ \Gr_{3}^{W}H^{3}(U,\bQ) $,
we use the spectral sequence
$$
E_{2}^{p,q} = H^{p}(S',R^{q}f_{*}\bQ_{U}) \Rightarrow
H^{p+q}(U,\bQ).
$$
Here we may assume that
$ S' $ is a non empty affine open subvariety of
$ S $ over which
$ f $ is smooth (shrinking
$ U $ and
$ S' $ if necessary) so that the
$ R^{q}f_{*}\bQ_{U} $ are variations of mixed Hodge structures.
The assumption on the Abel-Jacobi map for
$ S $ is equivalent to the vanishing of
$ \Gr_{3}^{W}H^{3}(S',\bQ) $ for
$ S' $ sufficiently small.
So it is enough to show
$$
H^{3-q}(S',\Gr_{q}^{W}R^{q}f_{*}\bQ_{U}) = 0\quad\text{for
$ 1 \le q \le 3 $,}
$$
using the spectral sequence associated to the filtration
$ W $, because
$ R^{q}f_{*}\bQ_{U} $ has weights
$ \ge q $ and
$ H^{i}(S',\Gr_{j}^{W}R^{q}f_{*}\bQ_{U}) $ has weights
$ \ge i + j $.

We have
$ H^{0}(S',\Gr_{3}^{W}R^{3}f_{*}\bQ_{U}) = 0 $ by the hypothesis
on the invariant part if general fibers have dimension
$ > 2 $, and the surface case follows from the dual of the
weak Lefschetz theorem (here we may assume
$ X_{s} $ projective).
The hypothesis on the
$ 2 $-forms implies
$ \Gr_{2}^{W}R^{2}f_{*}\bQ_{U} = 0 $ using the Lefschetz
theorem for divisors.
Finally
$ H^{2}(S',\Gr_{1}^{W}R^{1}f_{*}\bQ_{U}) = 0 $ by the
hypothesis in the case
$ \dim S > 1 $ (here the curve case is clear because
$ S' $ is affine).
So we get the assertion.

\medskip\noindent
{\bf 3.11.~Proof of Theorem (0.5).}
By the decomposition theorem [7], there is a noncanonical
isomorphism
$$
\bold{R}f_{*}\bQ_{X}[\dim X] \simeq \moplus_{i,Z} K_{Z}^{i}[-i],
\leqno(3.11.1)
$$
where
$ Z $ runs over the irreducible closed subvarieties of
$ S $ and the
$ K_{Z}^{i} $ are intersection complexes with local system
coefficients on
$ Z $, see loc. cit.
This decomposition holds in the derived category of mixed Hodge
Modules [47] on
$ S $, and
$ K_{Z}^{i} $ is pure of weight
$ i $, and is generically a variation of Hodge structure of weight
$ i - \dim Z $ (shifted by
$ \dim Z $).
Let
$ U $ be the open subvariety of
$ S $ over which
$ f $ is smooth.
Put
$ n = \dim X - \dim S $,
$ m = \dim S $.
Then
$ (K_{S}^{i}|_{U})[-m] $ is the local system
$ R^{i+n}f_{*}\bQ_{X}|_{U} $.

Let
$ \zeta \in \CH_{\AJ}^{2}(X)_{\bQ} $.
We have to show that this comes from
$ S $.
The restriction of
$ \zeta $ to a general fiber of
$ f $ is zero by hypothesis.
Using the localization sequence and spreading out [9],
there is a divisor
$ \Sigma $ on
$ S $ together with
$ \zeta' \in \CH^{1}(Z)_{\bQ} $ such that
$ i_{*}\zeta' = \zeta $ in
$ \CH^{2}(X)_{\bQ} $ where
$ Z = f^{-1}(\Sigma ) $ with the inclusion
$ i : Z \to X $.
So, by the injectivity of the cycle map in the divisor case (0.3),
it is enough to show that its cycle class in Deligne local
cohomology
$ H_{\cD,Z}^{4}(X,\bQ(2)) $ vanishes modifying
$ \zeta' $ by an element of
$ \CH^{1}(Z)_{\bQ} $ whose image in
$ \CH^{2}(X)_{\bQ} $ comes from
$ S $, because Deligne local cohomology is identified with
Deligne homology.
Here Deligne cohomology (or homology) means the absolute Hodge
cohomology (or homology), and we omit
$ {}^{\prime\prime} $ to simplify the notation.
This holds also for the later part of this section.

We define Deligne (local) cohomology with
coefficients in mixed Hodge Modules by
$$
\aligned
H_{\cD,\Sigma}^{j-i}(S,K_{S}^{i-n}(2)[-m])
&=\Hom_{D^{b}\MHM(S)}(\bQ_{S},i'_{*}i^{\prime !}
K_{S}^{i-n}(2)[j-i-m]),
\\
H_{\cD}^{j}(Z,K_{Z}^{i}(2))
&=\Hom_{D^{b}\MHM(Z)}(\bQ_{Z},K_{Z}^{i}(2)[j]),
\endaligned
$$
where
$ i' : \Sigma \to S $ denotes the inclusion.
We may assume that
$ \Sigma \supset S \setminus U $ replacing
$ \Sigma $ if necessary.
Then, by (3.11.1),
$ H_{\cD,Z}^{4}(X,\bQ(2)) $ is isomorphic to
$$
(\moplus_{i\le 2} H_{\cD,\Sigma}^{4-i}(S,K_{S}^{i-n}(2)[-m]))\oplus
(\moplus_{i,Z\ne S} H_{\cD}^{4-n-m-i}(Z,K_{Z}^{i}(2))).
\leqno(3.11.2)
$$
Here
$ H_{\cD,\Sigma}^{4-i}(S,K_{S}^{i-n}(2)[-m])) = 0 $ for
$ i > 2 $, because
$$
{}^{p}H^{j}i^{\prime !}K_{S}^{i-n} = 0 \quad\text{for
$ j \ne 1 $,}
\leqno(3.11.3)
$$
see [7] for
$ {}^{p}H^{j} $.
(Indeed, it vanishes for
$ j \ne 0,1 $ by the localization sequence, because the direct
image by an affine open embedding is an exact functor of
perverse sheaves.
The vanishing for
$ j = 0 $ follows from the property of intersection complex that
it has no nontrivial subobjects with strictly smaller support.)

Comparing (3.11.2) with a similar decomposition for
$ H_{\cD}^{4}(X,\bQ(2)) $,
we see that the cycle class of
$ \zeta' $ in
$ H_{\cD,Z}^{4}(X,\bQ(2)) $ is given by
$$
\msum_{i\le 2} \xi_{S}^{i} \in
\moplus_{i\le 2} H_{\cD,\Sigma}^{4-i}(S,K_{S}^{i-n}(2)[-m]),
$$
because
$ \zeta \in \CH_{\AJ}^{2}(X)_{\bQ} $.
Then the assertion follows from (0.3) if
$ \xi_{S}^{i} = 0 $ for any
$ i $.
Since
$ K_{S}^{-n}[-m] = \bQ_{S} $, we may assume
$ \xi_{S}^{0} = 0 $ modifying
$ \zeta' $ by a cycle coming from
$ S $ (using (0.3)).
Thus it is enough to show
$ \xi_{S}^{i} = 0 $ for
$ i = 1, 2 $, modifying
$ \zeta' $ by an element of
$ \CH^{1}(Z)_{\bQ} $ whose image in
$ \CH^{2}(X)_{\bQ} $ vanishes.

For
$ i = 2 $, the variation of Hodge structure
$ R^{2}f_{*}\bQ_{X}|_{U} $ has level
$ 0 $, and hence it is associated to an orthogonal representation
which has a finite monodromy group.
So we may assume that
$ K_{S}^{2-n}[-m] $ is a constant sheaf on a smooth variety
$ S $, replacing
$ X $ with a resolution of singularities of the base change of
$ X $ by a generically finite morphism to
$ S $, because the composition of the pull-back and the
pushforward of cycles under a generically finite morphism of
irreducible proper smooth varieties is the multiplication by
the generic degree.

Since
$ K_{S}^{2-n}[-m] $ is constant, there are cycles
$ \Gamma_{j} \in \CH^{1}(X)_{\bQ} $ such that the stalk
of
$ K_{S}^{2-n}(1)[-m] $ at
$ s \in U $ (i.e.
$ H^{2}(X_{s},\bQ)(1)) $ is generated by the cycle classes of the
restrictions of
$ \Gamma_{j} $ to
$ X_{s} $,
where
$ X_{s} = f^{-1}(s) $.
By an argument similar to [49, II]
(using the nearby cycle functor), we see that the cycle
class of
$ \Gamma_{j}|_{Z} $ in
$ H_{\cD,Z}^{4}(X,\bQ(2)) $ comes from the cycle class of
$ \Gamma_{j} $ in
$ H_{\cD}^{2}(X,\bQ(1)) $ using the composition of
$$
H_{\cD}^{2}(X,\bQ(1))\to H_{\cD}^{2}(Z,\bQ(1))\to
H_{\cD,Z}^{4}(X,\bQ(2)),
\leqno(3.11.4)
$$
where the last morphism is induced by the canonical morphism
$ \bQ_{Z} \to \bold{R}\Gamma_{Z}\bQ_{X}(1)[2] $.
Furthermore, using a decomposition similar to (3.11.2) for
$ H_{\cD}^{2}(X,\bQ) $,
we see that the image of (3.11.4) is contained in the direct
sum of
$ H_{\cD,\Sigma}^{4-i}(S,K_{S}^{i-n}(2)[-m]) $ with
$ i \le 2 $ under the decomposition (3.11.2).

Since
$ K_{S}^{2-n}(2)[-m] = \moplus_{j} \bQ_{S}(1) $ by using
$ \Gamma_{j} $, we have
$$
H_{\cD,\Sigma}^{2}(S,K_{S}^{2-n}(2)[-m]) =
\moplus_{j} H_{\cD,\Sigma}^{2}(S,\bQ(1)) =
\moplus_{j,k} H_{\cD,\Sigma_{k}}^{2}(S,\bQ(1)),
$$
where the
$ \Sigma_{k} $ are irreducible components of
$ \Sigma $ and
$ H_{\cD,\Sigma_{k}}^{2}(S,\bQ(1)) = \bQ $.
So the morphism to
$ H_{\cD,\Sigma}^{2}(S,K_{S}^{2-n}(2)[-m]) $ is identified with
the cycle class map of
$ [\Sigma_{k}] $ to
$ H_{\cD}^{2}(S,\bQ(1)) $.
If the cycle class of
$ \eta = \msum_{k} a_{k}[\Sigma_{k}] $ vanishes in
$ H_{\cD}^{2}(S,\bQ(1)) $,
we see that
$ \Gamma_{j} \cdot f^{*}\eta $ is rationally
equivalent to zero.
Applying this to each factor of
$ \xi_{S}^{2} $ we may assume
$ \xi_{S}^{2} = 0 $ by modifying
$ \zeta' $ modulo rational equivalence on
$ X $.

For
$ i = 1 $, we may assume
$ K_{S}^{1-n} \ne 0 $, because the assertion is clear otherwise.
Then the variation of Hodge structure
$ K_{S}^{1-n}[-1] $ is constant on
$ S \,(= \bP^{1}) $ by hypothesis.
Its stalk is given by
$ H = H^{1}(X_{s},\bQ) $ for a general
$ s \in S $, and is isomorphic to
$ H^{1}(X,\bQ) $.
Since
$ H_{\Sigma}^{2}(S,\bQ) = \moplus_{k}H_{\Sigma_{k}}^{2}(S,\bQ) =
\moplus_{k} \bQ(-1) $, we have
$$
\aligned
H_{\cD,\Sigma}^{3}(S,K_{S}^{1-n}(2)[-1])
&= \moplus_{k}\Ext_{\MHS}^{1}(\bQ,H(1)),
\\
H_{\cD}^{3}(S,K_{S}^{1-n}(2)[-1])
&= \Ext_{\MHS}^{1}(\bQ,H(1)),
\endaligned
$$
and
$ H_{\cD,\Sigma}^{3}(S,K_{S}^{1-n}(2)[-1]) \to
H_{\cD}^{3}(S,K_{S}^{1-n}(2)[-1]) $ is identified with the tensor
of the degree map
$ \moplus_{k} \bQ[\Sigma_{k}] \to \bQ $ with
$ \Ext_{\MHS}^{1}(\bQ,H(1)) \,(=\Pic(X)_{\bQ}^{0}) $.
So its kernel is generated by cycles rationally equivalent to
zero on
$ S \,(= \bP^{1}) $ tensored with
$ \Pic(X)_{\bQ}^{0} $,
and we may assume
$ \xi_{S}^{1} = 0 $ replacing
$ \zeta' $ modulo rational equivalence on
$ X $.
This completes the proof of (0.5).

\medskip\noindent
{\bf 3.12.~Theorem.}
{\it Let
$ X $ be an irreducible smooth proper complex algebraic variety.
Let
$ f : X \to S $ be a surjective morphism to
$ S := \bP^{1} $.
Assume that general fibers
$ X_{s} $ of
$ f $ are connected and the restriction morphism
$ H^{1}(X,\bQ) \to H^{1}(X_{s},\bQ) $ is an isomorphism for a
general
$ s \in S $.
Then the Albanese map for
$ X $ is injective if it is injective for
general fibers
$ X_{s} $.
}

\medskip\noindent
{\it Proof.}
The argument is similar to (3.11).
It is enough to consider the Albanese map tensored with
$ \bQ $.
Let
$ \zeta \in \CH_{\AJ}^{n}(X)_{\bQ} $ with
$ n = \dim X $.
We may assume that
$ f(\supp \zeta) $ is contained in the open subvariety
$ U $ over which
$ f $ is smooth.
Let
$ S' $ be a multivalued section of
$ f $ (i.e. finite over
$ S) $.
Let
$ d $ be its degree over
$ S $.
Then
$ f^{*}f_{*}\zeta $ and
$ [S'] \cdot f^{*}f_{*}\zeta $ are rationally equivalent
to zero, and we may assume that
$ f_{*}\zeta = 0 $ without any equivalence relation, modifying
$ \zeta $ by
$ d^{-1}[S'] \cdot f^{*}f_{*}\zeta $.

Now
$ \zeta $ determines an element
$ \xi_{s} $ of
$ \Alb(X_{s})\otimes_{\bZ}\bQ $ for each
$ s \in f(\supp \zeta) $.
By hypothesis the Albanese varieties
$ \Alb(X_{s}) $ form a constant abelian scheme over
$ U $ (using duality), and any section is constant because
$ S = \bP^{1} $.
Furthermore, its generic fiber is the Albanese variety of the
generic fiber
$ X_{K} $ of
$ f $ (where
$ K = \bC(S) $), and the Albanese map for
$ X_{K} $ tensored with
$ \bQ $ is surjective (using the base changes by finite
extensions of
$ K $).
So we can apply an argument similar to the case
$ i = 1 $ in (3.11), and we may assume that
$ \xi_{s} = 0 $ modifying
$ \zeta $.
Then
$ \zeta = 0 $ by hypothesis, and the assertion follows.

\medskip\noindent
{\bf 3.13.~Proof of Theorem (0.4).}
The injectivity of the maps follows by applying (0.5), (3.12)
inductively to a Lefschetz pencil.
The surjectivity of the Albanese map is clear.
For cycles of codimension
$ 2 $ it follows from (3.10) or [26] (using the
injectivity).

\medskip\noindent
{\bf 3.14.~Theorem.} {\it
Let
$ X $ be an irreducible smooth proper complex algebraic variety
with a surjective morphism
$ f : X \to S $ to a smooth curve
$ S $.
Assume that general fibers of
$ f $ are connected and have no global
$ 2 $-forms, and the Abel-Jacobi map for cycles of codimension
$ 2 $ on general fibers is injective.
Then Nori's conjecture {\rm [42]} holds, i.e.
$ \CH_{\AJ}^{2}(X) \subset \CH_{\alg}^{2}(X) $.
}

\medskip\noindent
{\it Proof.}
The argument is similar to (3.11) except that we consider the
cycle class in the usual cohomology instead of Deligne cohomology,
and rational equivalence is replaced by algebraic equivalence.
For
$ i = 1 $, we also use the vanishing of
$ H_{\Sigma}^{j}(S,K_{S}^{1-n}[-1]) $ for
$ j \ne 2 $, which follows from (3.11.3) because
$ \dim \Sigma = 0 $.
Then the assertion follows by an argument similar to (3.11).

\medskip\noindent
{\bf 3.15.~Proof of Theorem (0.6).}
The argument is similar to (3.11).
By the weak Lefschetz theorem,
$ \zeta $ belongs to
$ \CH_{\AJ}^{2}(X) $.
Since it is enough to show
$ \pi^{*}\zeta = 0 $ in
$ \CH^{2}(\tX) $, we may replace
$ X $ with
$ \tX $.
So
$ \tX $ will be denoted by
$ X $,
and
$ \pi^{*}\zeta $ by
$ \zeta $ from now on.

By spreading out there exist
$ \Sigma $ and
$ \zeta' $ as in (3.11) such that
$ i'_{*}\zeta' = \zeta $ in the notation of (3.11).
Since
$ \dim X_{s} \ge 3 $,
the monodromy of
$ R^{2}f_{*}\bQ_{X} $ is trivial, and
$ K_{S}^{2-n}[-1] $ is a constant variation of
Hodge structure.
Hence the restriction morphism
$ H^{2}(X,\bQ) \to H^{2}(X_{s},\bQ) $ is surjective and the
Picard number of
$ X_{s} $ is independent of
$ s \in U $.
So the assertion follows by the same argument as in (3.11).

\medskip\noindent
{\bf 3.16.~Theorem.} {\it
With the notation of {\rm (0.6),} assume
$ \dim X \ge 3 $.
Then
$ \zeta \in \CH^{2}(X,1)_{\bQ} $ is zero if its restriction to
$ X_{s} $ vanishes for any
$ s \in S' $.
In particular, the higher Abel-Jacobi map
$$
cl : \CH^{2}(X,1)_{\bQ} \to J(H^{2}(X,\bQ)(2))
\leqno(3.16.1)
$$
is injective if this map is injective for any
$ X_{s} \,(s \in S') $.
}

\medskip\noindent
{\it Proof.}
The argument is similar to (3.11) and (3.15).
It is sufficient to show the first assertion.
Since
$ R^{1}f_{*}\bQ_{X} $ is a constant local system by hypothesis,
it does not contribute to the Deligne local cohomology by a
weight argument.
Then, using the decomposition (3.11.1), the assertion is
reduced to the isomorphism between
$ \CH^{1}(Z,1) $ and the corresponding Deligne homology
[32, 3.1].

\medskip\noindent
{\bf 3.17.~Remarks.}
(i) The injectivity of the higher Abel-Jacobi map (3.16.1)
implies
$ H^{1}(X,\bQ) = 0 $,
see [51, 5.2].
It is conjectured that the converse is also true.

\medskip
(ii) If the assumptions of (3.10) are further satisfied in (3.14)
(e.g. if
$ \dim X = 3 $), then algebraic and homological equivalences
coincide for cycles of codimension
$ 2 $ on
$ X $ by (3.10) and (3.14).
(It is conjectured that this should hold assuming only the
nonexistence of global
$ 3 $-forms.)
In the case
$ \dim X = 3 $ and
$ \dim S = 1 $, we can prove the assertion also by using
Cor. 2.3 in [1].
Indeed, let
$ U $ be an open subvariety over which
$ f $ is smooth, and put
$ Y = f^{-1}(U) $.
Then by hypothesis and spreading out ([9], [15]), the diagonal cycle
$ [Y] $ in
$ Y\times_{U}Y $ is rationally equivalent to a cycle
$ \Gamma_{1} + \Gamma_{2} $ (with rational coefficients)
such that the
$ j $-th projection of
$ \Gamma_{j} $ to
$ Y $ is contained in a divisor for
$ j = 1, 2 $ (shrinking
$ U $ if necessary).
Using the embeddings
$ Y\times_{U}Y \to Y\times Y \to X\times X $,
a similar assertion holds for the diagonal
$ X $ of
$ X\times X $, and we can apply the theory of Barbieri-Viale on
balanced varieties in loc. cit.
The key point is that
$ \Griff^{2}(X) $ is a quotient of
$ H_{\Zar}^{0}(X,\cH^{3}(\bZ(2))) $ by [14] using the local-to-global
spectral sequence, and we have the action of the diagonal on
$ H_{\Zar}^{0}(X,\cH^{3}(\bZ(2))) $ which vanishes up to torsion.

\medskip
(iii) Nori's conjecture stated in the introduction is
equivalent to the one in [42] modulo Grothendieck's
generalized Hodge conjecture.
Indeed, the last conjecture implies that the abelian part
of the image of the Abel-Jacobi map coincides with the
abelian variety corresponding to the largest Hodge structure
of level
$ \le 1 $ contained in
$ H^{3}(X,\bQ) $.

\medskip
(iv) It is not easy to generalize the argument in (3.11)
to the case
$ p_{g}(X_{s}) \ne 0 $ even if we assume
$ \Gamma(X, \Omega_{X}^{3}) = 0 $, because the existence of the
transcendental part of
$ H^{2}(X_{s},\bQ) $ makes the situation completely different
(e.g. the Picard number of
$ X_{s} $ is not constant).

\medskip
(v) By a well-known conjecture of Beilinson [6] and Bloch [9],
$ \CH_{\AJ}^{2}(X)_{\bQ} $ should be determined by
$ H^{2}(X,\bQ) $ (more precisely, it should be expressed by
$ \Ext^{2}(\bQ,H^{2}(X,\bQ)(2)) $ in the conjectural
category of mixed motives).
Let
$ Y $ be a general surface in
$ X $ which is an intersection of general hyperplane
sections.
If the Hodge conjecture is true, there is a cycle
$ \Gamma \in \CH^{2}(Y\times X)_{\bQ} $ such that the
composition of the restriction morphism
$ H^{2}(X,\bQ) \to H^{2}(Y,\bQ) $ with
$ \Gamma_{*} : H^{2}(Y,\bQ) \to H^{2}(X,\bQ) $ is
the identity on
$ H^{2}(X,\bQ) $,
and Nori's conjecture can be reduced to the surface case
(where the conjecture is trivial) if the conjecture of
Beilinson and Bloch is true.
Note that the last conjecture can be replaced by their conjecture
on the injectivity of the Abel-Jacobi map for smooth projective
varieties over number fields, see [51, 0.4].

\medskip
(vi) Let
$ \pi : X \to Y $ be a
$ \bP^{1} $-bundle over a smooth projective variety
$ Y $, which has a section
$ D $.
If
$ D' $ is a sufficiently very ample divisor on
$ Y $, then
$ D + \pi^{*}D' $ is very ample on
$ X $.
Let
$ Z $ be a general hyperplane section of it.
Then
$ Z $ is a section of
$ \pi $ by calculating the intersection with
$ \pi^{-1}(y) $; in particular,
$ Y = Z $.

\bigskip\bigskip
\centerline{{\bf 4. Higher Abel-Jacobi map}}

\bigskip
\noindent
{\bf 4.1.~Currents.}
For a complex manifold
$ X $ of dimension
$ n $,
let
$ \cC^{\ssbull}(X) $ denote the complex of currents on
$ X $ which has the Hodge filtration
$ F $ as usual.
Here we normalize
$ \cC^{\ssbull}(X) $ so that
$ \cC^{j}(X) = 0 $ for
$ i > 0 $ or
$ i < -2n $.
It has a structure of double complex such that the Hodge filtration is
given by the first degree.
We have a natural morphism
$ \cE^{\ssbull}(X)(n)[2n] \to \cC^{\ssbull}(X) $,
where
$ \cE^{i}(X) $ denotes the vector space of
$ C^{\infty} $
$ i $-forms on
$ X $.
Let
$ S^{\ssbull}(X) $ denote the complex of locally finite
$ C^{\infty} $ chains on
$ X $ where
$ S^{-j}(X) $ consists of locally finite
$ j $-chains so that
$ H^{j}(S^{\ssbull}(X)) = H^{j+2n}(X,\bZ)(n) $.
There is a natural morphism
$ \iota : S^{\ssbull}(X) \to \cC^{\ssbull}(X) $.
The differential
$ d $ of
$ S^{\ssbull}(X) $ is defined in a compatible way with that of
$ \cC^{\ssbull}(X) $.
So it differs from the usual boundary map
$ \partial $ by the sign
$ -(-1)^{\deg} $ due to the Stokes theorem.
Note that the differential of a current
$ \Phi $ is defined by
$ (d\Phi)(\omega) + (-1)^{\deg \Phi}\Phi (d\omega) = 0 $
for
$ C^{\infty} $ forms
$ \omega $ with compact supports.
For a smooth complex algebraic variety
$ X $, we will denote
$ S^{\ssbull}(X^{\an}) $,
$ \cC^{\ssbull}(X^{\an}) $ by
$ S^{\ssbull}(X) $,
$ \cC^{\ssbull}(X) $ to simplify the notation.

Let
$ \oX $ be a smooth proper complex algebraic variety of dimension
$ n $, and
$ D $ a divisor on
$ \oX $ with normal crossings such that each irreducible component
$ D_{j} $ is smooth.
Put
$ X = \oX \setminus D $.
Let
$ \tD^{(j)} $ be the disjoint union of the intersections of
$ j $ irreducible components as in [22, II].
Then we have naturally a double complex
$$
\to F^{k}\cC^{\ssbull}(\tD^{(j+1)})
\to F^{k}\cC^{\ssbull}(\tD^{(j)})
\to \cdots
\to F^{k}\cC^{\ssbull}(\tD^{(0)})
\to 0
$$
by the dual construction of [22, III]
(using the push-down of currents instead of the pull-back of forms),
and the associated single complex will be denoted by
$ F^{k}\cC^{\ssbull}(\oX\anD) $.

We have the weight filtration
$ W $ on
$ \cC^{\ssbull}(\oX\anD) $ such that
$ \Gr_{j}^{W}\cC^{\ssbull}(\oX\anD) = \cC^{\ssbull}(\tD^{(j)})[j] $.
We define similarly
$ S^{\ssbull}(\oX\anD) $ with the filtration
$ W $ such that
$ \Gr_{j}^{W}S^{\ssbull}(\oX\anD) = S^{\ssbull}(\tD^{(j)})[j] $.
Then we get the polarizable mixed Hodge complex
$ K(\oX\anD) $ defined by
$$
(S^{\ssbull}(X),
(S^{\ssbull}(\oX\anD)_{\bQ}, W),
(\cC^{\ssbull}(\oX\anD); F, W);
S^{\ssbull}(X)_{\bQ},
(\cC^{\ssbull}(\oX\anD), W)),
$$
which calculates the Borel-Moore homology of
$ X $.
By (1.2.1) we have a canonical isomorphism (see also [32], [35]):
$$
H_{\cD}^{i}(X,\bZ(k)) =
H^{i-2n}(\Gamma(D''_{\cH}(K(\oX\anD)(k-n))).
\leqno(4.1.1)
$$

\medskip\noindent
{\bf 4.2.~Cycle class.}
With the above notation, let
$ \zeta = \msum_{j} (Z_{j},g_{j}) \in \CH^{p}(X,1) $ as in (2.2).
Put
$ d = n - p $.
Let
$ \gamma_{j} $ be the closure of the inverse image by
$ g_{j} $ of
$$
\{z \in \bC \,| \,\Re\, z > 0, \Im\, z = 0\} \subset \bP^{1}.
$$
Using a triangulation, it is viewed as a topological chain.
We give it an orientation so that
$ \partial \gamma_{j} = \div g_{j} $.
Then
$ \gamma := \msum_{j} \gamma_{j} $ is a topological cycle on
$ Z := \cup_{j} Z_{j} $, and it belongs to
$ S^{-2d-1}(X) $.
Let
$ \tZ_{j} \to Z_{j} $ be a resolution of singularities such
that the divisor of the pull-back
$ \tg_{j} $ of
$ g_{j} $ to
$ \tZ_{j} $ has normal crossings.
Let
$ \pi_{j} : \tZ_{j} \to X $ denote its composition with
the inclusion
$ i_{j} : Z_{j} \to X $.
Then we have the push-down of currents
$ (\pi_{j})_{*} : \cC^{\ssbull}(\tZ_{j}) \to \cC^{\ssbull}(X) $.

Let
$ \log_{\Hv} \tg_{j} $ denote a locally integrable function on
$ \tZ_{j} $ which is defined by choosing a branch of
$ \log \tg_{j} $ on
$ \tZ_{j} \setminus \tga_{j} $ where
$ \tga_{j} $ is the pull-back of
$ \gamma_{j} $ to
$ \tZ_{j} $.
(Hv stands for Heaviside.)
Then it is a current on
$ \tZ_{j} $, and we can verify
$$
\aligned
d(\log_{\Hv} \tg_{j})
&= {\tg}_{j}^{-1}d\tg_{j} -(2\pi i)\iota\tga_{j},
\\
d({\tg}_{j}^{-1}d\tg_{j})
&= (2\pi i)\iota(\div\tg_{j}).
\endaligned
\leqno(4.2.1)
$$
For example, we get the first equality by using the integration
on the inverse image of
$$
\{z \in \bC \,|\, \min\{|\arg z|, |z|, |z|^{-1} \} > \varepsilon \}
$$
for
$ \varepsilon \to 0 $.
Note that
$ {\tg}_{j}^{-1}d\tg_{j} $ is a form with locally integrable
coefficients on
$ \tZ_{j} $ and
$ \msum_{j} (\pi_{j})_{*}({\tg}_{j}^{-1}d\tg_{j}) $ is a closed current.
We define
$$
\aligned
(i_{j})_{*}(g_{j}^{-1}dg_{j})
&=(\pi_{j})_{*}({\tg}_{j}^{-1}d\tg_{j}),
\\
(i_{j})_{*}(\log_{\Hv} g_{j})
&=(\pi_{j})_{*}(\log_{\Hv} \tg_{j}),
\\
gdg
&=\msum_{j} (i_{j})_{*}(g_{j}^{-1}dg_{j})
\in F^{-d}\cC^{-2d-1}(\oX),
\\
\log_{\Hv} g
&=\msum_{j} (i_{j})_{*}(\log_{\Hv} g_{j})
\in \cC^{-2d-2}(\oX).
\endaligned
$$
These are independent of the choice of
$ \tZ_{j} $.

Let
$ \oZ_{j} $ be the closure of
$ Z_{j} $ in
$ \oX $.
Then
$ g_{j} $ is identified with a rational function
$ \og_{j} $ on
$ \oZ_{j} $, and we can define
$ \og^{-1}d\og =
\msum_{j}(\oi_{j})_{*}(\og_{j}^{-1}d\og_{j}) $, etc. similarly,
where
$ \oi_{j} : \oZ_{j} \to \oX $ is the inclusion
morphism.

Let
$ \div\og = \msum_{j} \div\og_{j} $.
This is supported on
$ D $, and there is a cycle
$ (\div\og)^{(1)} $ on
$ \tD^{(1)} $ such that its image in
$ \oX $ coincides with
$ \div\og $.
Taking a triangulation,
$ (\div\og)^{(1)} $ can be viewed as an element of
$ S^{-2d}(\tD^{(1)}) $.
So we get
$$
(g^{-1}dg)^{\wedge} :=
(\og^{-1}d\og, -(2\pi i)\iota (\div\og)^{(1)})
\in F^{-d}\cC^{-2d-1}(\oX\anD)
$$
such that
$ d(g^{-1}dg)^{\wedge} = 0 $.
Let
$ \oga $ be the closure of
$ \gamma $ in
$ \oX $.
Since
$ \partial \oga = \div\og $,
we get
$$
\gamma^{\wedge} := (\oga, -(\div\og)^{(1)}) \in
S^{-2d-1}(\oX\anD)
$$
such that
$ d\iota \gamma^{\wedge} = 0 $.
Then
$ d(\log_{\Hv} g) \in \cC^{-2d-1}(\oX\anD) $
coincides with the sum of
$ -(2\pi i)\iota\gamma^{\wedge} $ and
$ (g^{-1}dg)^{\wedge} $ in
$ \cC^{-2d-1}(\oX\anD) $ by (4.2.1).

\medskip\noindent
{\bf 4.3.~Theorem.} {\it
With the above notation,
$ cl(\zeta) \in H_{\cD}^{2p-1}(X,\bZ(p))'' $ corresponds
by the isomorphism
$ (4.1.1) $ to
$$
(-(2\pi i)^{-d}\gamma,
-(2\pi i)^{-d}\gamma^{\wedge},
-(2\pi i)^{-d-1}(g^{-1}dg)^{\wedge};
0, (2\pi i)^{-d-1}\log_{\Hv} g)
\leqno(4.3.1)
$$
where these elements belong to
$ S^{-2d-1}(X)(-d) $,
$ S^{-2d-1}(\oX\anD)_{\bQ}(-d) $,
$ F^{-d}\cC^{-2d-1}(\oX\anD) $,
$ S^{-2d-2}(X)_{\bQ}(-d) $ and
$ \cC^{-2d-2}(\oX\anD) $ respectively.
}

\medskip\noindent
{\it Proof.}
Since the class of
$ (\div\og)^{(1)} $ in
$ H^{2p-2}(\tD^{(1)},\bQ)(p-1) $ is a Hodge cycle, we see that
(4.3.1) belongs to
$ H^{-2d-1}(\Gamma(D''_{\cH}(K(\oX\anD)(-d))) $, see Remark (1.2).
Let
$ \oZ $ be the closure of
$ Z $ in $ \oX $, and
$$
\Sigma = \cup_{j} \supp \div\og_{j} \cup \Sing \oZ.
$$
Then the canonical morphism
$$
H_{\cD}^{2p-1}(X,\bZ(p))'' \to
H_{\cD}^{2p-1}(X \setminus \Sigma,\bZ(p))''
\leqno(4.3.2)
$$
is injective by the localization sequence.
So we may replace
$ X $ with
$ X \setminus \Sigma $, and assume that
$ Z $ is smooth, and hence irreducible.
Here we may assume also that the closure
$ \oZ $ of
$ Z $ in
$ \oX $ is a good smooth compactification (i.e.
$ \oZ \setminus Z $ is a divisor with normal crossings) by taking further
blowing-ups if necessary, and that every irreducible component of
$ \oZ \setminus Z $ is contained by only one irreducible component of
$ \oX \setminus X $.
Then the isomorphism (4.1.1) is compatible with
the push-forward by the closed embedding
$ Z \to X $, and the assertion is reduced to the case
$ X = Z $,
$ p = 1 $.

By (1.3) we have isomorphisms
$$
H_{\cD}^{1}(X,\bZ)'' = H_{\cD}^{1}(X,\bZ) =
H^{1}(\oX,C_{\oX\anD}^{\ssbull}\bZ(1)),
$$
and the cycle map is calculated by using the commutative diagram
$$
\CD
\Gamma (X,\bG_{m}) @>{\sim}>>
H^{1}(\oX,C_{\oX\anD}^{\ssbull}\bZ(1))
\\
@VVV @VVV
\\
\Gamma (X^{\an},\cO_{{X}^{\an}}^{*}) @>{\sim}>>
H^{1}(X,C_{X}^{\ssbull}\bZ(1)),
\endCD
\leqno(4.3.3)
$$
where the isomorphism on the bottom row is induced by the canonical
quasi-isomorphism
$$
\cO_{{X}^{\an}}^{*} = C(\bZ_{X^{\an}}(1) \to
\cO_{X^{\an}}).
$$
Here we may replace
$ H^{1}(X,C_{X}^{\ssbull}\bZ(1)) $ with the cohomology
of the single complex associated with
$$
S^{\ssbull}(X)(-d) \to
\cC^{\ssbull}(X)/F^{-d}\cC^{\ssbull}(X).
$$
Then, using (1.2.2), it is enough to show that the image of
$ g \in \Gamma (X^{\an}, \cO_{{X}^{\an}}^{*}) $ in
this cohomology is represented by
$$
(-(2\pi i)^{-d}\gamma,\quad (2\pi i)^{-d-1}\log_{\Hv} g),
$$
because the vertical morphisms of (4.3.3) are injective.

We can verify this assertion by using a Cech resolution together
with the delta functions supported on faces of a triangulation of
$ X^{\an} $ compatible with
$ \gamma $.
(See [28] for the notion of integral current.)
Indeed, let
$ \cU = \{U_{i}\}_{i\in \Lambda} $ be an open covering of
$ X^{\an} $ such that
$ U_{i} $ are simply connected.
We will denote by
$ C_{\cU}^{\ssbull}\cF $ the Cech complex of a sheaf
$ \cF $ associated with the covering
$ \cU $,
where
$$
C_{\cU}^{i}\cF = \moplus_{|I|=i+1} \Gamma (U_{I},\cF)\quad
\text{with
$ U_{I} = \cap_{i\in I} U_{i} $ for
$ I \subset \Lambda $}.
$$
For
$ g \in \Gamma (X^{\an},{O}_{{X}^{\an}}^{*}) $, we have the
corresponding element in the cohomology of the single complex
associated with
$$
C_{\cU}^{\ssbull}\bZ_{X^{\an}}(1) \to
C_{\cU}^{\ssbull}\cO_{X^{\an}}
$$
(where the first term has degree one), and it is given by
$$
(\{(\log(g|_{U_{i}}) - \log(g|_{U_{j}}))|_{U_{i,j}}\}_{i>j},
\{\log(g|_{U_{i}})\}_{i})\in
C_{\cU}^{1}\bZ_{X^{\an}}(1)\oplus
C_{\cU}^{0}\cO_{X^{\an}}.
$$
We define
$ C_{\cU}^{i}\cC^{j}(X) $,
$ C_{\cU}^{i}{S}^{j}(X) $ similarly.
Then
$$
\{(2\pi i)^{-d-1}((\log_{\Hv} g)|_{U_{i}}) - \log(g|_{U_{i}}))\}_{i} \in
C_{\cU}^{0}\cC^{-2d-2}(X)
$$
belongs to the image of
$ C_{\cU}^{0}S^{-2d-2}(X)(-d) $.
So we get the assertion, using the triple complex
$$
C_{\cU}^{\ssbull}S^{\ssbull}(X)(-d) \to
C_{\cU}^{\ssbull}(\cC^{\ssbull}(X)/
F^{-d}\cC^{\ssbull}(X)).
$$
This completes the proof of (4.3).

\medskip\noindent
{\bf 4.4.~Remark.}
The cycle map (2.4.1) induces
$$
\CH^{p}(X,1) \to \Hom_{\MHS}(\bZ,H^{2p-1}(X,\bZ)(p)),
\leqno(4.4.1)
$$
and (4.3) implies that the image of
$ \zeta $ by this morphism is represented by
$$
-(2\pi i)^{-d}\gamma \quad\text{and}\quad -(2\pi i)^{-d-1}
\msum_{j}(i_{j})_{*}(g_{j}^{-1}dg_{j}).
$$
Let
$ \CH_{\hom}^{p}(X,1) $ be the kernel of (4.4.1).
Then, in the notation of (1.1), the cycle map (2.4.1) induces the
higher Abel-Jacobi map
$$
\CH_{\hom}^{p}(X,1) \to J(H^{2p-2}(X,\bZ)(p))
= \Ext_{\MHS}^{1}(\bZ, H^{2p-2}(X,\bZ)(p)).
\leqno(4.4.2)
$$

Assume
$ H^{2p-1}(X,\bQ) = 0 $ or
$ X $ proper.
Then
$ \CH^{p}(X,1)/\CH_{\hom}^{p}(X,1) $ is finite because the target of
(4.4.1) is torsion.
So (4.4.2) induces the higher Abel-Jacobi map
$$
\CH^{p}(X,1)_{\bQ} \to J(H^{2p-2}(X,\bQ)(p)).
\leqno(4.4.3)
$$
By (4.3) this is expressed explicitly as follows.
For
$ \zeta \in \CH_{\hom}^{p}(X,1) $, there exist a
$ C^{\infty} $ chain
$ \Gamma $ on
$ X $ and
$ \Xi \in F^{-d}\cC^{-2d-1}(\oX\anD) $ such that
$$
\partial \Gamma = \gamma, \quad
d\Xi = (g^{-1}dg)^{\wedge}.
$$
By (4.3) and (1.2.2), the image of
$ \zeta $ under the higher Abel-Jacobi map (4.4.3) is
represented by the current
$$
\Phi_{\zeta} = (2\pi i)^{-d-1}
(\msum_{j} (i_{j})_{*}\log_{\Hv} g_{j} + (2\pi i)\iota\Gamma -
\Xi|_{X}).
\leqno(4.4.4)
$$
(Note that
$ d\iota\Gamma = - \iota\gamma $ by Stokes.)
If
$ X $ is proper, it is enough to consider
$ \Phi_{\zeta}(\omega) $ for
$ C^{\infty} $ forms
$ \omega $ which are direct sums of forms of type
$ (i,j) $ with
$ i \ge d + 1 $,
because the dual of
$ H^{2p-2}(X,\bC)/F^{p}H^{2p-2}(X,\bC) $ is
$ F^{d+1}H^{2d+2}(X,\bC) $.
Then
$ \Xi $ can be neglected, and we get the higher Abel-Jacobi map in
[5], [37] (see also [32]).

\medskip\noindent
{\bf 4.5.~Remark.}
With the notation of (1.1) and (2.2), the image of
$ \CH_{\dec}^{p-1}(X,1)_{\bQ} $ by the higher Abel-Jacobi map
(4.4.3) is contained in
$$
J(N^{p-1}H^{2p-2}(X,\bQ)(p)) \subset
\hbox{\rm Hdg}^{p-1}(X)_{\bQ}\otimes_{\bZ}\bC^{*},
$$
where
$ N^{p-1}H^{2p-2}(X,\bQ) $ is the
$ \bQ $-submodule generated by algebraic cycle classes, and
$ \hbox{\rm Hdg}^{p-1}(X)_{\bQ} $ is the group of Hodge cycles
with rational coefficients.
This can be reduced to the case
$ p = 1 $ by using resolutions of singularities, see e.g. [40].
The induced map
$$
\CH_{\ind}^{p}(X,1)_{\bQ} \to
J(H^{2p-2}(X,\bQ)(p))/
\hbox{\rm Hdg}^{p-1}(X)_{\bQ}\otimes_{\bZ}
\bC^{*}
\leqno(4.5.1)
$$
is called the reduced Abel-Jacobi map.

By [43], [46], the kernel of (4.5.1) for
$ p = 2 $ is isomorphic to the cokernel of
$$
K_{2}(\bC(X))_{\bQ} \to
\underset{U}\to\varinjlim\,\Hom_{\MHS}(\bQ, H^{2}(U,\bQ)(2))
\leqno(4.5.2)
$$
where the morphism is given by
$ d \log \wedge d \log $ at the level of integral logarithmic forms,
and the inductive limit is taken over nonempty open subvarieties
$ U $ of $ X $.

Indeed, for a divisor
$ Z $ on
$ X $, we have
$ H_{Z}^{3}(X,\bQ) = H_{2n-3}^{\BM}(Z,\bQ)(-n) $ and
$$
\CH_{\ind}^{1}(Z,1)_{\bQ}=
\Hom_{\MHS}(\bQ,H_{Z}^{3}(X,\bQ)(2))
$$
by (0.3) and (1.1.2).
Therefore, if
$ Z $ is sufficiently large, we get a short exact sequence
$$
\aligned
0 \to H^{2}(X,\bQ)/N^{1}H^{2}(X,\bQ)
&\to H^{2}(X \setminus Z,\bQ)
\\
\to \Ker
&(H_{Z}^{3}(X,\bQ) \to H^{3}(X,\bQ)) \to 0,
\endaligned
$$
and from its associated long exact sequence we can deduce
$$
\aligned
&\Hom_{\MHS}(\bQ, H^{2}(X \setminus Z, \bQ)(2))
\\
&= \Ker(\CH_{\ind}^{1}(Z,1)_\bQ \to
J((H^{2}(X,\bQ)/N^{1}H^{2}(X,\bQ))(2))),
\endaligned
$$
because
$ \Hom_{\MHS}(\bQ,H^{j}(X,\bQ)(2)) = 0 $ for
$ j = 2, 3 $.
Then it is enough to take the inductive limit over
$ Z $, and divide it by the image of
$ K_{2}(\bC(X))_{\bQ} $ under the tame symbol.

Note that the morphism
$ H^{2}(X \setminus Z,\bQ) \to H_{Z}^{3}(X,\bQ) $ is
given by taking the residue of logarithmic forms (at least on
$ Z_{\reg} $, see [22]), and the residue of
$ d \log f \wedge d \log g $ coincides with the differential of the
logarithm of the image of
$ \{f,g\} $ by the tame symbol up to sign.
The image of (4.5.1) is countable for
$ p \ge 2 $ by the rigidity argument of A. Beilinson [4] and
S. M\"uller-Stach [40], and (4.5.1) is not necessarily injective for
$ p \ge 3 $ (see [20]).
For
$ p = 2 $, it is conjectured that
$ \CH_{\ind}^{2}(X,1)_{\bQ} $ should be countable by C. Voisin [55],
and that (4.5.2) should be surjective by A. Beilinson [5].

It does not seem easy to prove the last conjecture by using [39].
Indeed, let
$ I $ be the target of (4.5.2) with
$ \bQ $ replaced by
$ \bZ $, and
$ I' $ be the image of
$ I'' := K_{2}(\bC(X)) $ in
$ I $.
Then
$ I''/m = I'/m = I/m $ for any positive integer
$ m $ by loc. cit. (using the exact sequence (3.7.1)).
Hence
$ I/I' $ is divisible.
It is torsion-free by the snake lemma, because so is
$ I $.
Therefore,
$ I/I' $ is uniquely divisible as proved in [46]
and we cannot get any more information.

\bigskip\bigskip
\centerline{{\bf 5. Construction of indecomposable higher cycles}}

\bigskip
\noindent
{\bf 5.1.~Elliptic surfaces.} Let
$ \pi : X \to C $ be an elliptic surface over a smooth proper
curve (i.e.,
$ X $ is smooth,
$ \pi $ is proper, and general fibers of
$ \pi $ are elliptic curves).
Let
$ \Sigma $ denote the smallest subset of
$ C $ such that
$ \pi $ is smooth over
$ C' := C \setminus \Sigma $.
Put
$ X' = \pi^{-1}(C') $,
$ X_{c} = \pi^{-1}(c) $ for
$ c \in C $.
It is well-known that the higher direct image sheaf
$ R^{1}\pi_{*}\bQ_{X} $ is an intersection complex (up to a shift)
by the decomposition theorem [7], and
$ H^{1}(C,R^{1}\pi_{*}\bQ_{X}) $ is a pure Hodge structure of
weight
$ 2 $,
see [57].
Let
$ L $ be the Leray filtration on
$ H^{2}(X,\bQ) $ which is defined by
$$
\aligned
L^{1}H^{2}(X,\bQ)
&= \cap_{c\in C} \Ker(H^{2}(X,\bQ) \to H^{2}(X_{c},\bQ)),
\\
L^{2}H^{2}(X,\bQ)
&= \pi^{*}H^{2}(C,\bQ),
\endaligned
$$
and
$ L^{0}H^{2}(X,\bQ) = H^{2}(X,\bQ), L^{3}H^{2}(X,\bQ)=0 $.
Then
$$
\Gr_{L}^{j}H^{2}(X,\bQ) = H^{j}(C,R^{2-j}\pi_{*}\bQ_{X}).
$$
We say that an elliptic surface has
{\it essentially no nontrivial section} if
$$
\Hom_{\MHS}(\bQ,H^{1}(C,R^{1}\pi_{*}\bQ_{X})(1)) = 0.
\leqno(5.1.1)
$$
By [57] this condition is equivalent to that the sections of
$ \pi : X \to C $ are torsion, if the monodromy of the local system
$ R^{1}\pi_{*}\bQ_{X}|_{C'} $ is nontrivial and
$ X $ is given a
$ 0 $-section.

Let
$ \{\pi_{t} : X_{t} \to C_{t}\}_{t\in \Delta} $ be a family of
elliptic surfaces over an open disk
$ \Delta $ of radius
$ \varepsilon $ whose restriction over
$ \Delta^{*}\,(:=\Delta\setminus \{0\}) $ is locally topologically
trivial.
We may view
$ \{\pi_{t} : X_{t} \to C_{t}\}_{t\in \Delta} $ as a small
deformation of
$ \pi_{0} : X_{0} \to C_{0} $ or a degeneration of elliptic
surfaces.
We say that it is {\it cohomologically nondegenerate} if
$ \{\Gr_{L}^{1}H^{2}(X_{t},\bQ)\}_{t\in \Delta} $ is a constant
local system on
$ \Delta $.
Note that
$ \dim \Gr_{L}^{1}H^{2}(X_{t},\bQ) $ is constant if and only if
so is
$ \dim \Gr_{L}^{0}H^{2}(X_{t},\bQ) $.
Using the decomposition theorem [7], the latter dimension
is given by
$ \msum_{c\in C_{t}} (n_{t,c}-1) + 1 $, where
$ n_{t,c} $ is the number of irreducible components of
$ X_{t,c} $.

Let
$ I $ denote the open interval
$ (0,\varepsilon ) $ contained in
$ \Delta^{*} := \Delta \setminus \{0\} $.
Let
$ \{\zeta_{t}\}_{t\in I} $ be an analytic family of higher cycles
with
$ \zeta_{t} \in \CH^{2}(X_{t},1) $.
We say that it {\it degenerates topologically} at
$ t = 0 $,
if there are
$ C^{\infty} $ families
$ \{\gamma_{t}\}_{t\in I}, \{\Gamma_{t}\}_{t\in I} $ as in (4.4)
which degenerate to a point as
$ t \to 0 $.

\medskip\noindent
{\bf 5.2.~Theorem.}
{\it Let
$ \{\pi_{t} : X_{t} \to C_{t}\}_{t\in \Delta} $ and
$ \{\zeta_{t}\}_{t\in I} $ be as above.
Assume that
$ \{X_{t}\} $ is cohomologically nondegenerate,
$ X_{0} $ has essentially no nontrivial section,
$ \{\zeta_{t}\} $ degenerates topologically at
$ t = 0 $, and
$ \int_{\Gamma_{t}} \omega_{t} \ne 0 $ for some
$ \omega_{t} \in \Gamma (X_{t},{\Omega}_{{X}_{t}}^{2}) $ where
$ t \in I $ is generic and
$ \Gamma_{t} $ is as above.
Then for a general
$ t \in I $,
the transcendental part of the image of
$ \zeta_{t} $ by the reduced Abel-Jacobi map {\rm (4.5.1)}
does not vanish {\rm (}i.e. its image in the Jacobian is not
contained in the image of
$ F^{1}H^{2}(X_{t},\bC)) $,
and hence
$ \zeta_{t} \ne 0 $ in
$ \CH_{\ind}^{2}(X_{t},1)_{\bQ} $.
}

\medskip\noindent
{\it Proof.}
Let
$ \cV_{t} = \Gr_{F}^{0}\Gr_{L}^{1}H^{2}(X_{t},\bQ) $.
Then
$ \{\cV_{t}\}_{t\in \Delta} $ is a holomorphic vector bundle
on
$ \Delta $.
By integrating holomorphic
$ 2 $-forms on
$ \Gamma_{t} $,
we get an analytic section
$ \sigma_{\Gamma} $ of
$ \{\cV_{t}\}_{t\in I} $,
because
$ H^{2}(C_{t},\bQ) $ and the image of
$ H^{2}(X_{t},\bQ) \to H^{2}(X_{t,c},\bQ) $ are of type
$ (1,1) $.
Let
$ L $ denote also the dual filtration on
$ H_{2}(X_{t},\bQ) $ such that
$ \Gr_{j}^{L}H_{2}(X_{t},\bQ) $ is the dual of
$ \Gr_{L}^{j}H^{2}(X_{t},\bQ) $,
i.e.
$$
\aligned
L_{0}H_{2}(X_{t},\bQ)
&= \Im(\moplus_{c} H_{2}(X_{t,c},\bQ)\to H_{2}(X_{t},\bQ)),
\\
L_{1}H_{2}(X_{t},\bQ)
&= \Ker(\pi_{*} : H_{2}(X_{t},\bQ) \to H_{2}(C_{t},\bQ)).
\endaligned
$$

Let
$ \{\eta_{t}\}\in\{\Gr_{1}^{L}H_{2}(X_{t},
\bQ)(-1)\}_{t\in\Delta} $
be a continuous family of topological cycles with rational
coefficients.
Then it determines a holomorphic section
$ \sigma_{\eta} $ of
$ \{\cV_{t}\}_{t\in \Delta} $ by integrating forms on it.
By (4.4) the image of
$ \zeta_{t} $ by (4.4.3) in
$ J(H^{2}(X_{t},\bQ)(2)) $ modulo the image of
$ F^{1}H^{2}(X_{t},\bC) $ is given by
$ \sigma_{\Gamma} $ modulo the image of
$ H_{2}(X_{t},\bQ)(-1) $.
If it vanishes for any
$ t \in I $,
then there exists a family
$ \{\eta_{t}\}_{t\in\Delta} $ as above such that
$ \sigma_{\Gamma} $ coincides with
$ \sigma_{\eta} $ over
$ I $,
because
$ I $ is an uncountable set.
By hypothesis, the limit of
$ \sigma_{\Gamma} $ for
$ t \to 0 $ is zero, and so is the value of
$ \sigma_{\eta} $ at the origin.
This implies that
$ \eta_{0} $ belongs to
$$
\Hom_{\MHS}(\bQ,\Gr_{L}^{1}H^{2}(X_{0},\bQ)(1))
$$
using Poincar\'e duality.
Therefore
$ \eta_{0} = 0 $ by (5.1.1), and hence
$ \eta_{t} = 0 $ for any
$ t \in \Delta $ by the triviality of the local system.
Thus we get
$ \sigma_{\Gamma} = 0 $.
But this contradicts the nonvanishing of
$ \int_{\Gamma_{t}} \omega_{t} $.

\medskip\noindent
{\bf 5.3.~Construction.}
Let
$ Y'_{t} = \bC^{2}, S' = \bC $,
and
$ g_{t} : Y'_{t} \to S' $ be a polynomial map defined by
$ g_{t} = y^{2} - x^{2}(x + t) $ for
$ t \in \bC $.
Taking the closure of the graph in
$ \bP^{2} \times S' $,
we get
$ \og_{t} : \overline{Y}'_{t} \to S' $,
and this gives an elliptic surface
$ f_{t} : Y_{t} \to S := \bP^{1} $ by taking the minimal
model of the singular fiber over
$ \infty \in \bP^{1} $ using Kodaira's classification.

If
$ t \ne 0 $,
the singular fibers over
$ 0 $ and
$ -4t^{3}/27 $ are both rational curves with one ordinary double
point.
(For example, if we put
$ s = y/x $,
then
$ x = s^{2} - t, y = s(s^{2} - t) $ on
$ \{g_{t}(x,y) = 0\} $.)
So each of these singular fibers determines a higher cycle by taking
the normalization and choosing a rational function with simple zero
and pole at the pull-back of the double point.
This is well-defined up to a sign and modulo decomposable cycles.
It
$ t $ is positive, we can choose the rational function on the
singular fiber at
$ 0 $ so that
$ \gamma $ and
$ \Gamma $ in (4.2) and (4.4) are given respectively by
$$
\aligned
\gamma'_{t}
&= \{(x,y) \in \bR^{2} : g_{t}(x,y) = 0, x \le 0\},
\\
\Gamma'_{t}
&= \{(x,y) \in \bR^{2} : g_{t}(x,y) \le 0, x \le 0\},
\endaligned
$$
where the function is
$ -(y-\sqrt{t}x)/(y+\sqrt{t}x) = -(s-\sqrt{t})/(s+\sqrt{t}) $.

The higher cycles constructed above
are still decomposable because
$ Y_{t} $ is a rational surface, and we have to take a base
change.
Let
$ \rho : C \to S $ be a generic hyperelliptic curve such that
$ \infty \in S \,(= \bP^{1}) $ is a ramification point, but
$ 0 \in S $ is not.
We assume that
$ H^{1}(C,\bQ) $ does not contain a Hodge structure isomorphic
to the cohomology of the elliptic curve defined by
$ y^{2} = x^{3} + 1 $ (with
$ j $-invariant
$ 0) $.
Let
$ \tf_{t} : X_{t}\to C $ be the minimal model of the base change of
$ f_{t} : Y_{t}\to S $ by
$ \rho $.
(Actually, we can also consider the open surface with the singular
fiber over
$ \infty $ deleted, because
$ 1 $ is not an eigenvalue of the local monodromy of
$ R^{1}(f_{t})_{*}\bQ_{Y_{t}}|_{S'} $ around
$ \infty $ so that
$ R^{1}(f_{t})_{*}\bQ_{Y_{t}} =
\bold{R}j_{*}j^{-1}R^{1}(f_{t})_{*}\bQ_{Y_{t}} $ where
$ j : S' \to S $ denotes the inclusion,
and similarly for the pull-back by
$ \rho $.)
We assume that
$ \rho $ is not ramified over
$ 0 $, and choose a point
$ \tzero $ of
$ C $ over
$ 0 $.
Let
$ \gamma_{t} $ be the connected component of the pull-back of
$ \gamma'_{t} $ contained in the fiber over
$ \tzero $.
There is a connected component
$ \Gamma_{t} $ of the pull-back of
$ \Gamma'_{t} $ such that
$ \partial\Gamma_{t} = \gamma_{t} $ for
$ t $ sufficiently small.

\medskip\noindent
{\bf 5.4~Remarks.}
(i) The image of
$ \Gamma'_{t} $ by
$ g_{t} $ is the interval
$ [-4t^{3}/27,0] $, and
$ \Gamma'_{t} $ gives a degeneration of
$ \gamma'_{t} \subset Y_{t,0} $ as
$ c \to -4t^{3}/27 $, i.e.
$ \gamma'_{t} $ is the vanishing cycle associated to the
singular fiber over
$ -4t^{3}/27 $.
The cohomological nondegeneration is related to the phenomenon
that two
$ A_{1} $-singularities appear by a deformation of a holomorphic
function with an isolated singularity of type
$ A_{2} $, i.e. two
$ A_{1} $-singularities of a function can join and degenerate
to an
$ A_{2} $-singularity.

\medskip
(ii) Instead of the hyperelliptic curve
$ C $,
it is also possible to consider an
$ m $-fold cyclic covering of
$ \bP^{1} $ ramified over two points
$ \alpha $ and
$ \infty $ if
$ \alpha $ is generic and
$ m $ is prime to
$ 6 $ and strictly greater than
$ 6 $.

\medskip\noindent
{\bf 5.5.~Theorem.}
{\it If
$ \varepsilon $ is sufficiently small, then
$ \{X_{t}\} $ satisfies the assumptions of {\rm (5.2).}
}

\medskip\noindent
Indeed, the hypotheses are satisfied by the following lemmas
and proposition:

\medskip\noindent
{\bf 5.6.~Lemma.}
{\it The family
$ \{X_{t}\} $ is cohomologically nondegenerate.
}

\medskip\noindent
{\it Proof.}
For
$ t \in \Delta $ we have
$$
H^{1}(S',R^{1}(f_{t})_{*}\bQ_{Y_{t}}|_{S'}) =
H^{1}(S',R^{1}(g_{t})_{*}\bQ_{Y'_{t}}) = 0,
$$
see e.g. [24].
In particular, these groups are constant for
$ t \in \Delta $,
and this holds also for the cohomology of its restriction over
a small open disk with center
$ 0 $ in
$ S' $ if
$ \varepsilon $ is sufficiently small, because there are no
singular fibers over the complement of the disk in
$ S' $ if
$ t $ is sufficiently small.
(Note that by hypothesis, the direct image sheaf is defined over
$ \cup S'_{t} $ and its restriction to each
$ S'_{t} $ is
$ R^{1}(g_{t})_{*}\bQ_{Y'_{t}} $;
moreover, a similar assertion holds for the pull-back by
$ \rho $.)
So we can use the Mayer-Vietoris sequence to show that the
cohomology of the pull-back of the sheaf to
$ C $ is constant, and the assertion follows.

\medskip\noindent
{\bf 5.7.~Lemma.}
{\it The elliptic surface
$ X_{0} $ has essentially no nontrivial section.
}

\medskip\noindent
{\it Proof.}
The local system associated to
$ f_{0} $ has a finite monodromy group, and is trivialized by
taking the pull-back under a finite base change.
This holds also for the local system associated to
$ X_{0} $,
and it is enough to show that there is essentially no nontrivial
section (see (5.1.1)) for the constant elliptic surface over
$ C $ whose fiber
$ X_{c} $ is defined by the equation
$ y^{2} = x^{3} + 1 $.
But the self-duality of
$ H^{1}(X_{c},\bQ) $ implies
$$
\Hom(\bQ,H^{1}(X_{c},\bQ)\otimes H^{1}(C,\bQ)(1)) =
\Hom(H^{1}(X_{c},\bQ),H^{1}(C,\bQ)),
$$
and the assertion follows from the hypothesis.

\medskip\noindent
{\bf 5.8.~Proposition.}
{\it The integral
$ \int_{\Gamma_{t}} \omega_{t} $ does not vanish for some holomorphic
$ 2 $-from
$ \omega_{t} $.
}

\medskip\noindent
{\it Proof.}
We may assume that
$ C $ is given by the equation
$ s^{2} = h(z) $ with
$$
h(z) = \prod_{j=1}^{2g+1}(z - \alpha_{j}),
$$
where
$ z $ is the coordinate of
$ S' $, the
$ \alpha_{j} $ are generic complex numbers, and
$ g \ge 1 $.
Then
$ dz/\sqrt{h(z)} $ is a nonzero
$ 1 $-form on
$ C $.
We show that there is a nonzero
$ 2 $-form
$ \omega_{t} $ on
$ X_{t} $ whose restriction to
$ Y'_{t}\times_{S'}C' $ is
$ dx\wedge dy/f_{t}^{*}\sqrt{h(z)} $,
where
$ C' $ is the inverse image of
$ S' $.
(This is rather trivial if we assume that the genus of
$ C $ is sufficiently large, compared with the order of the pole of
$ dx\wedge dy/f_{t}^{*}dz $.)

Let
$ \omega_{\rel} = dx\wedge dy/f_{t}^{*}dz $.
It gives a section of
$ (f_{t})_{*}\omega_{Y_{t}/S} $ by [36], because
$ \omega_{\rel} $ satisfies the Gauss hypergeometric differential
equation
$$
z(z + 4t^{3}/27)\partial_{z}^{2}\omega_{\rel}+
(2z + 4t^{3}/27)\partial_{z}\omega_{\rel}+(5/36)\omega_{\rel} = 0,
$$
(via the Gauss-Manin connection, see e.g. [24]), and the roots of
its indicial equation [18] at
$ 0, -4t^{3}/27, \infty $ are respectively
$ \{0,0\}, \{0,0\} $ and
$ \{1/6,5/6\} $.
Similarly,
$ \rho^{*}\omega_{\rel} $ is extended to a section of (the direct
image of)
$ \omega_{X_{t}/C} $,
because the roots of the indicial equation at
$ \rho^{-1}(\infty ) $ are
$ 1/3, 5/3 $.
Then it gives a nonzero section of
$ \omega_{X_{t}} $ by using the section
$ dz/\sqrt{h(z)} $ of
$ \omega_{C} $,
and it coincides with the above
$ 2 $-form
$ \omega_{t} $.

To show
$ \int_{\Gamma_{t}} \omega_{t} \ne 0 $,
we may assume that the
$ \alpha_{j} $ are real and sufficiently small so that
$ f_{t}^{*}\sqrt{h(z)} > 0 $ on
$ \Gamma_{t} $.
Then the assertion is clear.

\newpage
\centerline{{\bf References}}
\bigskip

\item{[1]}
L. Barbieri-Viale, Balanced varieties, in Algebraic
$ K $-theory and its applications (H. Bass et al. eds.),
World Sci. Publishing, River Edge, NJ, 1999, pp. 298--312.

\item{[2]}
L. Barbieri-Viale and V. Srinivas, A reformulation of Bloch's conjecture,
C.R. Acad. Sci. Paris 321 (1995), 211--214.

\item{[3]}
R. Barlow, Rational equivalence of zero cycles for some more surfaces
with
$ p_{g} = 0 $, Inv. Math. 79 (1985), 303--308.

\item{[4]}
A. Beilinson, Higher regulators and values of
$ L $-functions, J. Soviet Math. 30 (1985), 2036--2070.

\item{[5]}
\SameAuthor, Notes on absolute Hodge cohomology, Contemporary Math. 55
(1986) 35--68.

\item{[6]}
\SameAuthor, Height pairing between algebraic cycles, Lect. Notes in
Math., vol. 1289, Springer, Berlin, 1987, pp. 1--26.

\item{[7]}
A. Beilinson, J. Bernstein and P. Deligne, Faisceaux pervers,
Ast\'erisque, vol. 100, Soc. Math. France, Paris, 1982.

\item{[8]}
S. Bloch, Torsion algebraic cycles and a theorem of Roitman,
Compos. Math. 39 (1979), 107--127.

\item{[9]}
\SameAuthor, Lectures on algebraic cycles, Duke University
Mathematical series 4, Durham, 1980.

\item{[10]}
\SameAuthor, Algebraic cycles and higher
$ K $-theory, Advances in Math., 61 (1986), 267--304.

\item{[11]}
\SameAuthor, Algebraic cycles and the Beilinson conjectures,
Contemporary Math. 58 (1) (1986), 65--79.

\item{[12]}
\SameAuthor, The moving lemma for higher Chow groups,
J. Alg. Geom. 3 (1994), 537--568.

\item{[13]}
S. Bloch, A. Kas and D. Lieberman, Zero cycles on surfaces with
$ p_{g} = 0 $, Compos. Math. 33 (1976), 135--145.

\item{[14]}
S. Bloch and A. Ogus, Gersten's conjecture and the homology of
schemes, Ann. Sci. Ecole Norm. Sup. 7 (1974), 181--201.

\item{[15]}
S. Bloch and V. Srinivas,
Remarks on correspondences and algebraic cycles,
Amer. J. Math. 105 (1983), 1235--1253.

\item{[16]}
J. Carlson, Extensions of mixed Hodge structures, in Journ\'ees
de G\'eom\'etrie Alg\'ebrique d'Angers 1979,
Sijthoff-Noordhoff Alphen a/d Rijn, 1980, pp. 107--128.

\item{[17]}
H. Clemens and P. Griffiths, The intermediate Jacobian of the cubic
threefold, Ann. of Math. 95 (1972), 281--356.

\item{[18]}
E. Coddington and N. Levinson, Theory of Ordinary Differential
Equations, McGraw-Hill, 1955.

\item{[19]}
A. Collino, Griffiths' infinitesimal invariant and higher
$ K $-theory on hyperelliptic Jacobians, J. Alg. Geom. 6 (1997),
393--415.

\item{[20]}
A. Collino and N. Fakhruddin, Indecomposable higher Chow cycles on
Jacobians, preprint.

\item{[21]}
P. del Angel and S. M\"uller-Stach,
The transcendental part of the regulator map for
$ K_{1} $ on a mirror family of
$ K3 $ surfaces, preprint.

\item{[22]}
P. Deligne, Th\'eorie de Hodge I, Actes Congr\`es Intern. Math., 1970,
vol. 1, 425--430; II, Publ. Math. IHES, 40 (1971), 5--57; III ibid., 44
(1974), 5--77.

\item{[23]}
C. Deninger and A. Scholl, The Beilinson conjectures, in Proceedings
Cambridge Math. Soc. (eds. Coats and Taylor) 153 (1992), 173--209.

\item{[24]}
A. Dimca and M. Saito,
Algebraic Gauss-Manin systems and Brieskorn modules, Am. J. Math. 123
(2001), 163--184.

\item{[25]}
F. El Zein and S. Zucker,
Extendability of normal functions associated to algebraic cycles,
in Topics in transcendental algebraic geometry, Ann. Math. Stud., 106,
Princeton Univ. Press, Princeton, N.J., 1984, pp. 269--288.

\item{[26]}
H. Esnault and M. Levine, Surjectivity of cycle maps,
Ast\'erisque 218 (1993), 203--226.

\item{[27]}
H. Esnault and E. Viehweg, Deligne-Beilinson cohomology, in Beilinson's
conjectures on Special Values of
$ L $-functions, Academic Press, Boston, 1988, pp. 43--92.

\item{[28]}
H. Federer, Geometric Measure Theory,
Springer, New York, 1969.

\item{[29]}
H. Gillet, Deligne homology and Abel-Jacobi maps,
Bull. Amer. Math. Soc. 10 (1984), 285--288.

\item{[30]}
P. Griffiths, On the period of certain rational integrals I,
II, Ann. Math. 90 (1969), 460--541.

\item{[31]}
H. Inose and M. Mizukami, Rational equivalence of
$ 0 $-cycles on some surfaces of general type with
$ p_{g} = 0 $. Math. Ann. 244 (1979),205--217.

\item{[32]}
U. Jannsen, Deligne homology, Hodge-$ D $-conjecture, and motives,
in Beilinson's conjectures on Special Values of
$ L $-functions, Academic Press, Boston, 1988, pp. 305--372.

\item{[33]}
\SameAuthor, Mixed motives and algebraic
$ K $-theory, Lect. Notes in Math., vol. 1400, Springer, Berlin, 1990.

\item{[34]}
\SameAuthor, Letter from Jannsen to Gross on higher Abel-Jacobi maps,
in Proceedings of the NATO Advanced Study Institute on The arithmetic
and geometry of algebraic cycles (B.B. Gordon et al. eds.), Kluwer
Academic, Dordrecht, 2000, pp. 261--275.

\item{[35]}
J.R. King, Log complexes of currents and functorial properties
of the Abel-Jacobi map, Duke Math. J. 50 (1983), 1--53.

\item{[36]}
J. Koll\'ar, Higher direct images of dualizing sheaves, I, II, Ann.
of Math. 123 (1986), 11--42; 124 (1986), 171--202.

\item{[37]}
M. Levine, Localization on singular varieties, Inv. Math. 91 (1988),
423--464.

\item{[38]}
\SameAuthor, Bloch's higher Chow groups revisited,
Ast\'erisque 226 (1994), 235--320.

\item{[39]}
A.S. Merkur'ev and A.A. Suslin,
$ K $-cohomology of Severi-Brauer varieties and the norm residue
homomorphism, Math. USSR Izv. 21 (1983), 307--340.

\item{[40]}
S. M\"uller-Stach, Constructing indecomposable motivic cohomology
classes on algebraic surfaces, J. Alg. Geom. 6 (1997), 513--543.

\item{[41]}
D. Mumford, Rational equivalence of
$ 0 $-cycles on surfaces, J. Math. Kyoto Univ. 9 (1969), 195--204.

\item{[42]}
M.V. Nori, Algebraic cycles and Hodge theoretic connectivity,
Inv. Math. 111 (1993), 349--373.

\item{[43]}
C. Pedrini, Bloch's conjecture and the $ K $-theory of projective
surfaces, in The arithmetic and geometry of algebraic cycles, CRM Proc.
Lecture Notes, 24, Amer. Math. Soc., Providence, 2000, pp. 195--213.

\item{[44]}
A. Roitman, Rational equivalence of zero cycles, Math. USSR Sbornik 18
(1972), 571--588.

\item{[45]}
\SameAuthor, The torsion in the group of zero cycles modulo rational
equivalence, Ann. Math. 111 (1980), 553--569.

\item{[46]}
A. Rosenschon, Indecomposable Elements of
$ K_{1} $, K-theory 16 (1999), 185--199.

\item{[47]}
M. Saito, Mixed Hodge Modules, Publ. RIMS, Kyoto Univ., 26 (1990),
221--333.

\item{[48]}
\SameAuthor, On the formalism of mixed sheaves, RIMS-preprint 784, Aug.
1991.

\item{[49]}
\SameAuthor, Hodge conjecture and mixed motives, I, Proc. Symp. Pure
Math. 53 (1991), 283--303; II, in Lect. Notes in Math., vol. 1479,
Springer, Berlin, 1991, pp. 196--215.

\item{[50]}
\SameAuthor, Mixed Hodge complex on algebraic varieties, Math. Ann.
316 (2000), 283--331.

\item{[51]}
\SameAuthor, Arithmetic mixed sheaves, Inv. Math. 144 (2001),
533--569.

\item{[52]}
R.I. Shafarevich, Algebraic surfaces, Proc. the Steklov Institute
of Mathematics, No. 75 (1965), translated by American Mathematical
Society, Providence, R.I.

\item{[53]}
C. Schoen, Some examples of torsion in the Griffiths group,
Math. Ann. 293 (1992), 651--679.

\item{[54]}
C. Voisin, Sur les z\'ero cycles de certaines hypersurfaces munies
d'un automorphisme, Ann. Sci. Norm. Sup. Pisa 19 (4) (1992),
473--492.

\item{[55]}
\SameAuthor, Remarks on zero-cycles of self-products of varieties,
in Moduli of Vector Bundles, Lect. Notes in Pure and Applied
Mathematics, vol. 179, M. Dekker, New York, 1996, pp. 265--285.

\item{[56]}
\SameAuthor, The Griffiths group of a general Calabi-Yau threefold
is not finitely generated, Duke Math. J. 102 (2000), 151--186.

\item{[57]}
S. Zucker, Hodge theory with degenerating coefficients,
$ L_{2} $-cohomology in the Poincar\'e metric, Ann. Math.,
109 (1979), 415--476.

\bigskip\noindent
\ver
\bye